\documentclass[10 pt]{article}
\usepackage{amsmath,amsfonts,amssymb,amsthm}
\usepackage[dvips]{graphics}
\usepackage[english]{babel}
\usepackage{makeidx}
\usepackage{latexsym}
\usepackage[latin1]{inputenc}
\newtheorem{thm}{Theorem}[section]
\newtheorem{cor}{Corollary}[section]
\newtheorem{lem}{Lemma}[section]
\newtheorem{prop}{Proposition}[section]
\newcommand{\<}{\langle}
\renewcommand{\>}{\rangle}

\newcommand{\N}[1][]{\ensuremath{{\mathbb{N}^{#1}} }}
\newcommand{\Z}[1][]{\ensuremath{{\mathbb{Z}^{#1}} }}
\newcommand{\R}[1][]{\ensuremath{{\mathbb{R}^{#1}} }}
\newcommand{\T}[1][]{\ensuremath{{\mathbb{T}^{#1}} }}
\newcommand{\Q}[1][]{\ensuremath{{\mathbb{Q}^{#1}} }}

\begin{document}
\baselineskip=1.0\baselineskip
\title{\textbf{Global Well-Posedness and Non-linear Stability of Periodic Traveling Waves
for a Schr\"odinger-Benjamin-Ono System}}
\author{Jaime Angulo P. $^1$, Carlos Matheus $^2$, and  Didier Pilod $^3$ }
\footnotetext[1]{Partially supported by CNPq/Brazil; email:
angulo@ime.unicamp.br.} \footnotetext[2]{Partially supported by
CNPq/Brazil; email: matheus@impa.br} \footnotetext[3]{Partially
supported by CNPq/Brazil; email: pilod@impa.br}

\date{August 1, 2007}

\maketitle \large{\vspace{-0.5cm}

{\scriptsize \centerline{$^1$Department of Mathematics, IME-USP}
 \centerline{Rua do Matão 1010, Cidade Universit\'aria, CEP 05508-090, São Paulo, SP, Brazil.}
 \centerline{$^2$ IMPA, Estrada Dona Castorina 110,}
\centerline{ CEP 22460-320 Rio de Janeiro, RJ,
 Brazil.}
\centerline{$^3$UFRJ, Institute of Mathematics, Federal University
of Rio de Janeiro,} \centerline{P.O. Box 68530 - Cidade
Universit\'aria - Ilha do Fund\~ao.} \centerline{CEP 21945-970 Rio
de Janeiro, RJ, Brazil.}}

\vspace{0.5cm}

\begin{abstract} The objective of this paper is two-fold: firstly,
we develop a local and global (in time) well-posedness theory for a
system describing the motion of two fluids with different densities
under capillary-gravity waves in a deep water flow (namely, a
Schr\"odinger-Benjamin-Ono system) for \emph{low-regularity} initial
data in both periodic and continuous cases; secondly, a family of
new periodic traveling waves for the Schr\"odinger-Benjamin-Ono
system is given:  by fixing a minimal period we obtain, via the
implicit function theorem, a smooth branch of periodic solutions
bifurcating a Jacobian elliptic function called {\it dnoidal}, and,
moreover, we prove that all these periodic traveling waves are
nonlinearly stable by perturbations with the same wavelength.
\end{abstract}


\begin{center}
\section{Introduction}
\end{center}
In this paper we are interested in the study of the following
Schr\"odinger-Benjamin-Ono (SBO) system
\begin{equation}\label{eq:lsws}
\left \{
    \begin{array}{l}
iu_t+u_{xx}=\alpha vu,\\
v_t+\gamma Dv_x=\beta (|u|^2)_x,
 \end{array}
   \right.
\end{equation}
where $u$ is a complex-valued function, $v$ is a real-valued
function,  $t \in \Bbb R$, $x \in \Bbb R$ or $\Bbb T$, $\alpha,
\beta$ and $\gamma$ are real constants such that $\alpha \neq 0$ and
$\beta \neq 0$, and $D\partial_x$ is a linear differential operator
representing the dispersive term. Here $D=\mathcal H\partial_x$
where $\mathcal H$ denotes the {\it Hilbert transform} defined as
$$
\widehat{\mathcal H f}(k)=-i\text{sgn}(k)\widehat{f}(k),
$$
where \begin{equation*} \text{sgn}(k)  =\left\{
  \begin{array}
  [c]{r}
  -1,\text{ }k<0,\\
  1,\text{ }k>0.
  \end{array}
\right.
\end{equation*}
Note that from these definitions we have that $D$ is a linear
positive Fourier operator with symbol $|k|$. The system
(\ref{eq:lsws}) was deduced by Funakoshi and Oikawa (\cite{fo:fo}).
It describes the motion of two fluids with different densities under
capillary-gravity waves in a deep water flow. The short surface wave
is usually described by a Schr\"odinger type equation and the long
internal wave is described by some sort of wave equation accompanied
with a dispersive term (which is a Benjamin-Ono type equation in
this case). This system is also of interest in the sonic-Langmuir
wave interaction in plasma physics \cite{k:k}, in the
capillary-gravity interaction wave \cite{dr:dr}, \cite{g:g}, and in
the general theory of water wave interaction in a nonlinear medium
\cite{bn1:bn1}, \cite{bn2:bn2}. We note that the Hilbert transform
considered in \cite{fo:fo} for describing system (\ref{eq:lsws}) is
given as $-\mathcal H$.

When studying an initial value problem, the first step is usually to
investigate in which function space well-posedness occurs. In our
case, smooth solutions of the SBO system (\ref{eq:lsws}) enjoy the
following conserved quantities
\begin{equation}\label{eq:law2}
\left \{
    \begin{array}{l}
G(u,v)\equiv Im\,\int\,u(x)
\overline{u_x(x)}\,dx+\frac{\alpha}{2\beta}\int\,|v(x)|^2\,dx,\\
\\
E(u,v)\equiv\int\,|u_x(x)|^2\,dx +\alpha
\int\,v(x)|u(x)|^2\,dx-\frac{\alpha\gamma}{2\beta}
\int\,|D^{1/2}v(x)|^2\;dx,\\
\\
H(u,v)\equiv\int\,|u(x)|^2\,dx,
\end{array}
\right.
\end{equation}
where $D^{1/2}$ is the Fourier multiplier defined as $\widehat
{D^{1/2}v}(k)=|k|^{\frac12}\widehat v(k)$. Therefore, the natural
spaces to study well-posedness are the Sobolev $H^s$-type spaces.
Moreover, due to the scaling property of the SBO system
(\ref{eq:lsws}) (see \cite{bop} Remark 2), we are led to investigate
well-posedness in the spaces $H^s \times H^{s-1/2}$, $s \in \Bbb R$.

In the  continuous case Bekiranov, Ogawa and Ponce \cite{bop1}
proved local well-posedness for initial data in $H^s(\Bbb R)\times
H^{s-\frac12}(\Bbb R)$ when $|\gamma| \neq 1$ and $s\geq 0$. Thus,
because of the conservation laws in (\ref{eq:law2}), the solutions
extend globally in time when $s\geq 1$, in the case
$\frac{\alpha\gamma}{\beta}<0$. Recently, Pecher \cite{per} has
shown local well-posedness in $H^s(\Bbb R)\times H^{s-\frac12}(\Bbb
R)$ when $|\gamma| = 1$ and $s> 0$. He also used the Fourier
restriction norm method to extend the global well-posedness result
when $1/3<s<1$, always in the case $\frac{\alpha\gamma}{\beta}<0$.
Here, we improve the global well-posedness result till $L^2(\mathbb
R)\times H^{-\frac12}(\mathbb R)$ in the case $\gamma \neq 0$ and
$|\gamma| \neq 1$. Indeed, we refine the bilinear estimates of
Bekiranov, Ogawa and Ponce \cite{bop1} in Bourgain spaces
$X^{0,b_1}\times Y_{\gamma}^{-\frac12,b}$ with $b, \ b_1<\frac12$
(see Proposition \ref{prop-g1}). These estimates combined with the
$L^2$-conservation law allow us to show that the size of the time
interval provided by the local well-posedness theory depends only on
the $L^2$-norm of $u_0$. It is worth to point out that this scheme
applies for other dispersive systems. In fact Colliander, Holmer and
Tzirakis \cite{cht} already applied this method to Zakharov and
Klein-Gordon-Schr\"odinger systems. They also announced the above
result for the SBO system (see Remark 1.5 in \cite{cht}). However
they allowed us to include it in this paper since there were not
planing to write it up anymore. Note that we also prove global
well-posedness in $H^s(\mathbb R)\times H^{s-\frac12}(\mathbb R)$
when $s>0$ in the case $\gamma \neq 0$ and $|\gamma|\neq1$. We take
the opportunity to express our gratitude to Colliander, Holmer and
Tzirakis for the fruitful interaction about the
Schr\"odinger-Benjamin-Ono system.

In the periodic setting, there does not exist, as far as we know,
any result about the well-posedness of the SBO system
(\ref{eq:lsws}). Bourgain \cite{Bo} proved well-posedness for the
cubic nonlinear Schr\"odinger equation (NLS) (see (\ref{sch})) in
$H^s(\Bbb T)$ for $s \ge 0$ using the Fourier transform restriction
method. Unfortunately, this method does not apply directly for the
Benjamin-Ono equation. Nevertheless, using an appropriate Gauge
transformation introduced by Tao \cite{Tao}, Molinet \cite{Mo}
proved well-posedness in $L^2(\Bbb T)$. Here we apply Bourgain's
method for the SBO system to prove its local well-posedness in
$H^s(\Bbb T)\times H^{s-\frac12}(\Bbb T)$ when $s \ge 1/2$ in the
case $\gamma\neq 0$, $|\gamma| \neq 1$. The main tool is the new
bilinear estimate stated in Proposition \ref{prop1}. Furthermore, by
standard arguments based on the conservation laws, this leads to
global well-posedness in the energy space $H^1(\T) \times
H^{1/2}(\T)$ in the case $\frac{\alpha\gamma}{\beta}<0$. We also
show that our results are sharp in the sense that the bilinear
estimates on these Bourgain spaces fail whenever $s< 1/2$ and
$|\gamma| \neq 1$ or $s \in \Bbb R$ and $|\gamma| =1$. In fact, we
use Dirichlet's Theorem on rational approximation to locate certain
plane waves whose nonlinear interactions behave badly in low
regularity.

In the second part of this paper, we turn our attention to another
important aspect of dispersive nonlinear evolution equations: the
traveling-waves. These solutions imply a balance between the effects
of nonlinearity and dispersion. Depending on the specific boundary
conditions on the wave's shape, these special states of motion can
arise as either solitary or periodic waves. The study of this
special steady waveform is essential for the explanation of many
wave phenomena observed in the practice: in surface water waves
propagating in a canal, in propagation of internal waves or in the
interaction between long waves and short waves as in our case. In
particular, some questions such as existence and stability of these
traveling waves are very important in the understanding of the
dynamic of the equation under investigation.

The solitary waves are in general a single crested, symmetric,
localized traveling waves, with sech-profiles (see Ono \cite{ono}
and Benjamin \cite{benjamin} for the existence of solitary waves of
algebraic type or with a finite number of oscillations). The study
of the nonlinear stability or instability of solitary waves has had
a big development and refinement in recent years. The proofs have
been simplified and sufficient conditions have been obtained to
insure the stability to small localized perturbations in the
waveform. Those conditions have showed to be effective in a variety
of circumstances, see for example \cite{albert1}, \cite{albert2},
\cite{albert3}, \cite{b:b}, \cite{bona1}, \cite{gss:gss},
\cite{w1:w1}.

The situation regarding to the study of periodic traveling waves is
very different. The stability and the existence of explicit formulas
of these progressive wave trains have received comparatively little
attention. Recently many research papers about this issue have
appeared for specify dispersive equations, such as  the existence
and  stability of {\it cnoidal waves} for the Korteweg-de Vries
equation \cite{ABS} and the stability of {\it dnoidal waves} for the
one-dimensional cubic nonlinear Schr\"odinger equation
\begin{equation}\label{sch}
iu_t+u_{xx}+|u|^2u=0,
\end{equation}
where $u=u(t,x)\in{\Bbb C}$ and  $x,t\in{\Bbb R}$ (Angulo
\cite{a:a}, see also Angulo\&Linares \cite{AL} and
Gallay\&H$\breve{a}$r$\breve{a}$gus \cite{gahara1}, \cite{gahara2}).

In this paper we are also interested in giving a stability theory of
periodic traveling waves solutions for the nonlinear dispersive
system SBO (\ref{eq:lsws}). The periodic traveling waves solutions
considered here will be of the general form
\begin{equation}\label{eq:ptw}
\left \{
    \begin{array}{l}
u(x,t)=e^{i\omega t}e^{ic(x-ct)/2}\phi(x-ct),\\
v(x,t)=\psi(x-ct),
\end{array}
   \right.
\end{equation}
where $\phi,\psi:\Bbb R\to \Bbb R$ are smooth, $L$-periodic
functions (with a prescribed period $L$), $c>0$, $\omega\in \Bbb R$
and we will suppose that there is a $q\in\Bbb N$ such that
$$
4q\pi/{c}=L.
$$
So, by replacing these permanent waves form  into (\ref{eq:lsws}) we
obtain  the pseudo-differential system
\begin{equation}\label{eq:eds}
\left \{
    \begin{array}{l}
\phi''-\sigma\phi=\alpha \,\psi\phi\\
\gamma\mathcal H\psi'-c\psi=\beta \phi^2 +A_{\phi,\psi}
\end{array}
   \right.
   \end{equation}
where $\sigma=\omega-\frac{c^2}{4}$ and $A_{\phi,\psi}$ is an
integration constant which we will set equal zero in our theory.
Existence of analytic solutions of system (\ref{eq:eds}) for
$\gamma\neq 0$ is a difficult task. In the framework of traveling
waves of type solitary waves, namely, the profiles  $\phi, \psi$
satisfy $\phi(\xi), \psi(\xi)\to 0$ as $|\xi|\to \infty$, it is well
known the existence of solutions for (\ref{eq:eds}) in  the form
\begin{equation}\label{eq:sol}
 \phi_{0,s}(\xi)=\sqrt{\frac{2c\sigma}{\alpha\beta}}sech(\sqrt
\sigma \xi),\qquad \psi_{0,s}(\xi)=-\frac{\beta}{c}\phi_{0,s}^2(\xi)
\end{equation}
when  $\gamma=0$, $\sigma>0$, and $\alpha\beta>0$. For $\gamma\neq
0$ a theory of even solutions of these permanent waves solutions has
been established in \cite{am1:am1} (see also \cite{am2:am2}) by
using the concentration-compactness method.

For $\gamma=0$ and $\sigma>2\pi^2/{L^2}$ we prove (along the lines
of Angulo \cite{a:a} with regard to (\ref{sch})) the existence of a
smooth curve of even periodic traveling wave solutions for
(\ref{eq:eds}) with $\alpha=1,\beta=1/2$; note that this restriction
does not imply loss of generality. This construction is based on the
{\it dnoidal} Jacobian elliptic function , namely,
\begin{equation}\label{eq:dno}
\left \{
    \begin{array}{l}
\phi_0(\xi)=\eta_1\;dn\Big(\frac{\eta_1}{2\sqrt{c}}\;\xi;k\Big)\\
\psi_0(\xi)=-\frac{\eta_1^2}{2c}\;dn^2\Big(\frac{\eta_1}{2\sqrt{c}}\;\xi;k\Big),
\end{array}
\right.
\end{equation}
where $\eta_1$ and $k$ are positive smooth functions depending of
the parameter $\sigma$. We observe that the solution in
(\ref{eq:dno}) gives us in `` the limit " the solitary waves
solutions (\ref{eq:sol}) when $\eta_1\to \sqrt{4c\sigma}$ and $k\to
1^-$, because in this case the elliptic function $dn$ converges,
uniformly on compacts sets, to the hyperbolic function $sech$.

In the case of our main interest, $\gamma\neq 0$, the existence of
periodic solutions is a delicate issue. Our approach for the
existence of these solutions uses the implicit function theorem
together with the explicit formulas in (\ref{eq:dno}) and a detailed
study of the periodic eigenvalue problem associated to the Jacobian
form of Lame's equation
\begin{equation}\label{eq:2.21}
\left \{
    \begin{array}{l}
 \frac{d^2}{dx^2}\Psi+[\rho-6k^2sn^2(x;k)]\Psi=0\\
\Psi(0)=\Psi(2K(k)),\;\;\Psi'(0)=\Psi'(2K(k)),
\end{array}
\right.
\end{equation}
where  $sn(\cdot;k)$ is the Jacobi elliptic function of type
snoidal and $K=K(k)$ represents the complete elliptic integral of
the first kind and defined for $k\in (0,1)$ as
$$
K(k)=\int_0^1\;\frac{dt}{\sqrt{(1-t^2)(1-k^2t^2)}}.
$$
So, by fixing a period $L$, and choosing $c$ and $\omega$ such that
$ \sigma\equiv \omega- \frac{c^2}{4}$ satisfies $\sigma
>2\pi^2/{L^2}$, we obtain a smooth branch
$\gamma\in (-\delta, \delta)\to (\phi_\gamma, \psi_\gamma)$ of
periodic traveling wave solutions of (\ref{eq:eds}) with a
fundamental period $L$ and bifurcating from $(\phi_0, \psi_0)$ in
(\ref{eq:dno}). Moreover, we obtain that for $\gamma$ near zero
$\phi_\gamma(x)>0$ for all $x\in \Bbb R$ and $\psi_\gamma(x)<0$ for
$\gamma<0$ and $x\in \Bbb R$.

Furthermore, concerning the non-linear stability of this branch of
periodic solutions, we extend the classical approach developed by
Benjamin \cite{b:b}, Bona \cite{bona1} and Weinstein \cite{w1:w1} to
the periodic case. In particular, using the conservation laws
(\ref{eq:law2}), we prove that the solutions $(\phi_\gamma,
\psi_\gamma)$ are stable in $H^1_{per}([0,L])\times
H^{\frac12}_{per}([0,L])$ at least when $\gamma$ is negative near
zero. We use essentially the Benjamin\&Bona\&Weinstein's stability
ideas because it gives us an easy form of manipulating with the
required spectral conditions and the positivity property of the
quantity $\frac{d}{d\sigma}\int\phi^2_\gamma(x)dx$, which are basic
information in our stability analysis.

However, we do not use the abstract stability theory of Grillakis
{\it et al.} in our approach basically because of the two
circumstances above. We recall that Grillakis {\it et al.} theory in
general requires a study of the Hessian for the function
$$
d(c,\omega)= L (e^{ic\xi/2}\phi_\gamma,\psi_\gamma)\equiv E
(e^{ic\xi/2}\phi_\gamma,\psi_\gamma)+c G
(e^{ic\xi/2}\phi_\gamma,\psi_\gamma)+\omega H
(e^{ic\xi/2}\phi_\gamma,\psi_\gamma)
$$
with $\gamma=\gamma(c,\omega)$, and a specific spectrum information
of the matrix linear operator $H_{c,\omega}= L''
(e^{ic\xi/2}\phi_\gamma,\psi_\gamma)$. In our case, these facts  do
not seem to be easily obtained.

So, for $\gamma<0$ we reduce the required spectral information (see
formula~(\ref{eq:3.6})) to the study of the self-adjoint operator
$\mathcal L_\gamma$,
$$
\mathcal L_\gamma= -\frac{d^2}{d\xi^2} +\sigma
+\alpha\psi_\gamma-2\alpha\beta \phi_\gamma\circ\mathcal
K_\gamma^{-1}\circ\phi_\gamma,
$$
where $\mathcal K_\gamma^{-1}$ is the inverse operator of $\mathcal
K_\gamma= -\gamma D+c$. Hence we obtain via the min-max principle
that $\mathcal L_\gamma$ has a simple negative eigenvalue and zero
is a simple eigenvalue with eigenfunction $\frac{d}{dx}\phi_\gamma$
provide that $\gamma$ is small enough.


Finally, we close this introduction with the organization of this
paper: in Section 2, we introduce some notations to be used
throughout the whole article; in Section 3, we prove the global
well-posedness results in the periodic and continuous settings via
some appropriate bilinear estimates; in Section 4, we show the
existence of periodic traveling waves by the implicit function
theorem; then, in Section 5, we derive the stability of these waves
based on the ideas of Benjamin and Weinstein, that is, to manipulate
the information from the spectral theory of certain self-adjoint
operators and the positivity of some relevant quantities.

\begin{center}
\item\section{Notation}
\end{center}
\setcounter{section}{2} \setcounter{equation}{0}

For any positive numbers $a$ and $b$, the notation $a \lesssim b$
means that there exists a positive constant $\theta$ such that $a
\le \theta b$. Here, $\theta$ may depend only on certain parameters
related to the equation (\ref{eq:lsws}) such as $\gamma$, $\alpha$,
$\beta$. Also, we denote $a \sim b$ when, $a \lesssim b$ and $b
\lesssim a$.

For $a \in \mathbb{R}$, we denote by $a+$ and $a-$ a number slightly
larger and smaller than $a$, respectively.

In the sequel, we fix $\psi$ a smooth function supported on the
interval $[-2,2]$ such that $\psi(x)\equiv 1$ for all $|x|\leq 1$
and, for each $T>0$, $\psi_T(t):=\psi(t/T)$.

Let $L>0$, the inner product of two functions in $L^2([0,L]) $ is
given by
\begin{displaymath}
<f,g>=\int_0^Lf(x)\bar g(x)dx, \quad \forall f, \ g \in L^2([0,L]).
\end{displaymath}
Now let $\mathcal P'_L$ the set of periodic distributions of period
$L$, for all $s \in \mathbb R$ we denote by
$H^s_{per}([0,L])=H^s_L(\mathbb R)$ the set of all $f$ in $\mathcal
P'_L$ such that
\begin{displaymath}
\|f\|_{H^s_L}=\left(L\sum_{n=-\infty}^{+\infty}(1+|n|^2)^{s}
|\widehat{f}(n)|^2\right)^{\frac12}<\infty,
\end{displaymath}
where $(\widehat{f}(n))_{n \in \mathbb Z}$ denote the Fourier series
of $f$ (for further information see Iorio\&Iorio \cite{ii:ii}).
Sometimes we also write $H^s(\mathbb{T})$ to denote the space
$H^s_{per}([0,L])$ when the period $L$ does not play a fundamental
role.

Similarly, when $s \in \mathbb R$, we denote by $H^s(\mathbb R)$ the
set of all $f \in \mathcal{S}'(\mathbb R)$ such that
\begin{displaymath}
\|f\|_{H^s}=\left(\int_{\mathbb R}(1+|\xi|^2)^{s}
|\widehat{f}(\xi)|^2d\xi\right)^{\frac12}<\infty,
\end{displaymath}
where $\mathcal{S}'(\mathbb R)$ is the set of tempered distributions
and $\widehat{f}$ is the Fourier transform of $f$.

When the function $u$ is of the two time-space variables $(t,x) \in
\Bbb R \times \mathbb R$, periodic in space of period $L$, we define
its Fourier transform by
$$\widehat{u}(\tau,n)=\frac{1}{(2\pi)^{1/2}L}
\int_{\Bbb R \times [0,L]}u(t,x)e^{-i(n\pi x/\ell+t\tau)}dtdx,$$ and
similarly, when $u :\Bbb R\times \Bbb R \rightarrow \Bbb C$, we
define
$$\widehat{u}(\tau,\xi)=\frac{1}{2\pi}
\int_{\Bbb R^2}u(t,x)e^{-i(x\xi+t\tau)}dtdx.$$

Next, we introduce the Bourgain spaces related to the
Schr\"odinger-Benjamin-Ono system in the periodic case:
\begin{equation} \label{X}
\|u\|_{X_{per}^{s,b}} := \left( \sum_{n \in
\mathbb{Z}}\int_{\mathbb{R}}\langle\tau+n^2\rangle^{2b}\langle
n\rangle^{2s}|\widehat{u}(\tau, n)|^2d\tau \right)^{1/2},
\end{equation}
\begin{equation} \label{Y}
\|u\|_{Y^{s,b}_{\gamma,per}} := \left( \sum_{n \in
\mathbb{Z}}\int_{\mathbb{R}}\langle\tau+\gamma|n|n\rangle^{2b}\langle
n\rangle^{2s}|\widehat{u}(\tau, n)|^2d\tau \right)^{1/2},
\end{equation}
and the continuous case:
\begin{equation} \label{X'}
\|u\|_{X^{s,b}} := \left(
\int_{\mathbb{R}}\int_{\mathbb{R}}\langle\tau+\xi^2\rangle^{2b}\langle
\xi\rangle^{2s}|\widehat{u}(\tau, \xi)|^2 d\xi d\tau \right)^{1/2},
\end{equation}
\begin{equation} \label{Y'}
\|u\|_{Y^{s,b}_{\gamma}} := \left(
\int_{\mathbb{R}}\int_{\mathbb{R}}\langle\tau+\gamma|\xi|\xi\rangle^{2b}\langle
\xi\rangle^{2s}|\widehat{u}(\tau, \xi)|^2 d\xi d\tau \right)^{1/2},
\end{equation}
where $\langle x\rangle:=1+|x|$. The relevance of these spaces are
related to the fact that they are well-adapted to the linear part of
the system and, after some time-localization, the coupling terms of
(\ref{eq:lsws}) verifies particularly nice bilinear estimates.
Consequently, it will be a standard matter to conclude our global
well-posedness results (via Picard fixed point method).

\begin{center}
\item\section{Global Well-Posedness of the Schr\"odinger-Benjamin-Ono System}
\end{center}

This section is devoted to the proof of our well-posedness results
for (\ref{eq:lsws}) in both continuous and periodic settings.

\setcounter{section}{3} \setcounter{equation}{0} \indent
\begin{center}
{3.1 \ \it Global well-posedness on $\mathbb{R}$}
\end{center}

The bulk of this subsection is the proof of the following theorem:

\begin{thm} \label{thm-g1}
Let $0<|\gamma|\neq 1$. Then, the SBO system is globally well-posed
for initial data $(u_0,v_0) \in H^s(\mathbb R)\times
H^{s-\frac12}(\mathbb R)$, when $s \ge 0$.
\end{thm}

In the rest of this section, we will denote by
$U(t):=e^{it\partial_x^2}$ and $V_{\gamma}(t):=e^{-\gamma t\mathcal
H \partial_x^2}$ the unitary groups associated to the linear part of
\eqref{eq:lsws}. The proof of Theorem \ref{thm-g1} follows the lines
of \cite{cht}. Let us first state the linear estimates:
\begin{lem} \label{lem-g3}
Let $0\le b, b_1\le\frac12$, $s \in \mathbb R$ and $0 < T  \le 1$.
Then
\begin{equation} \label{lem-g3.1}
\|\psi_TU(t)u_0\|_{X^{s,b_1}} \lesssim T^{\frac12-b_1}\|u_0\|_{H^s},
\end{equation}
and
\begin{equation} \label{lem-g3.2}
\|\psi_TV_{\gamma}(t)v_0\|_{Y^{s-\frac12,b}_{\gamma}} \lesssim
T^{\frac12-b}\|v_0\|_{H^{s-\frac12}},
\end{equation}
for $\gamma \in \mathbb R$.
\end{lem}

\noindent \textbf{Proof of Lemma \ref{lem-g3}.} Estimate
\eqref{lem-g3.1} is proved in  \cite{cht} Lemma 2.1 (a). Next we
combine Estimate \eqref{lem-g3.1} and the fact that
\begin{equation} \label{lem-g3.3}
V_{\gamma}(t)=P_{+}U(\gamma t)-P_{-}U(-\gamma t)
\end{equation}
to deduce Estimate \eqref{lem-g3.2}, where
\begin{displaymath}
\widehat{P_{+}f}=\chi_{(0,+\infty)}\widehat{f} \quad \text{and}
\quad \widehat{P_{-}f}=\chi_{(-\infty,0)}\widehat{f}.
\end{displaymath}
\hfill $\square$

\begin{lem} \label{lem-g4}
\begin{itemize}
\item[(i)]
Let $s \in \mathbb R$, $0 < T \le 1$, $0\le  c_1\le\frac12$ and $b_1
\ge 0$ such that  $b_1+c_1\le 1$. Then
\begin{equation} \label{lem-g4.1}
\|\psi_T\int_0^tU(t-t')z(t')dt'\|_{X^{s,b_1}} \lesssim
T^{1-b_1-c_1}\|z\|_{X^{s,-c_1}},
\end{equation}
and
\begin{equation} \label{lem-g4.2}
\|\int_0^tU(t-t')z(t')dt'\|_{C([0,T];H^s)} \lesssim
T^{\frac12-c_1}\|z\|_{X^{s,-c_1}}.
\end{equation}
\item[(ii)]
Let $s \in \mathbb R$, $\gamma \in \mathbb R$, $0 < T \le 1$, $0\le
c\le\frac12$ and $b \ge 0$ such that  $b+c\le 1$. Then
\begin{equation} \label{lem-g4.3}
\|\psi_T\int_0^tV_{\gamma}(t-t')z(t')dt'\|_{Y^{s-\frac12,b}_{\gamma}}
\lesssim T^{1-b-c}\|z\|_{Y^{s-\frac12,-c}_{\gamma}},
\end{equation}
and
\begin{equation} \label{lem-g4.4}
\|\int_0^tV_{\gamma}(t-t')z(t')dt'\|_{C([0,T];H^{s-\frac12})}
\lesssim T^{\frac12-c}\|z\|_{Y^{s-\frac12,-c}_{\gamma}}.
\end{equation}
\end{itemize}
\end{lem}

\noindent \textbf{Proof of Lemma \ref{lem-g4}.} For the proof of
Estimates \eqref{lem-g4.1} and \eqref{lem-g4.2} see \cite{cht} Lemma
2.3 (a). For Estimates \eqref{lem-g4.3} and \eqref{lem-g4.4} we
combine Identity \eqref{lem-g3.3} with Estimates \eqref{lem-g4.1}
and \eqref{lem-g4.2}. \hfill $\square$

Once these linear estimates are established, our task is to prove
the following bilinear estimates:

\begin{prop} \label{prop-g1}
Let $\gamma \in \mathbb R$ such that $|\gamma|\neq 1$ and $\gamma
\neq 0$. Then, we have for any $\frac14<b, \ b_1, \ c, \ c_1
<\frac12$
\begin{equation} \label{prop-g1.1}
\|uv\|_{X^{0,-c_1}} \lesssim
\|u\|_{Y^{-\frac12,b}_{\gamma}}\|v\|_{X^{0,b_1}}, \quad \text{if}
\quad b+b_1+c_1 \ge 1,
\end{equation}
\begin{equation} \label{prop-g1.2}
\|\partial_x(u\bar{v})\|_{Y^{-\frac12,-c}_{\gamma}} \lesssim
\|u\|_{X^{0,b_1}}\|v\|_{X^{0,b_1}}, \quad \text{if} \quad 2b_1+c \ge
1,
\end{equation}
where the implicit constants depend on $\gamma$.
\end{prop}
For the proof of these bilinear estimates, we need the following
standard Bourgain-Strichartz estimates:
\begin{prop} \label{prop-g2}
Let $\gamma \in \mathbb R$ such that $|\gamma|\neq 1$ and $\gamma
\neq 0$. Then
\begin{equation} \label{prop-g2.1}
\|u\|_{L^3_{t,x}} \lesssim \|u\|_{X^{0,1/4+}},
\end{equation}
and
\begin{equation} \label{prop-g2.2}
\|u\|_{L^3_{t,x}} \lesssim \|u\|_{Y^{0,1/4+}_{\gamma}}.
\end{equation}
\end{prop}

\noindent Finally, we recall the two following technical lemmas
proved in \cite{GTV}:
\begin{lem} \label{lem-g1} Let $f \in L^q(\mathbb R), \ g \in L^{q'}(\mathbb
R)$ with $1 \le q, \ q' \le +\infty$ and $\frac1q+\frac1q'=1$.
Assume that $f$ and $g$ are nonnegative, even and nonincreasing for
positive argument. Then $f\ast g$ enjoys the same property. In
particular $f\ast g$ takes its maximum at zero.
\end{lem}
\begin{lem} \label{lem-g2}
Let $0\le a_{1}, a_{2}<\frac12$ such that $a_1+a_2>\frac12$. Then
$$\int_{\mathbb R}\langle y-\alpha \rangle^{-2a_1}\langle y-\beta \rangle^{-2a_2}
\lesssim \langle \alpha -\beta\rangle^{1-2(a_1+a_2)}, \quad \forall
\; \alpha, \beta \in \mathbb R.$$
\end{lem}

After these preliminaries, we are ready to show the bilinear
estimates (\ref{prop-g1.1}) and (\ref{prop-g1.2}):

\medskip

\noindent \textbf{Proof of Proposition \ref{prop-g1}.} Without loss
of generality we can suppose that $|\gamma|<1$ in the rest of the
proof.

We first begin with the proof of Estimate \eqref{prop-g1.1}. Letting
$f(\tau, \xi)=\langle \xi\rangle^{-1/2}\langle\tau+\gamma
\xi|\xi|\rangle^{b}\widehat{u}(\tau, \xi)$, $g(\tau,
\xi)=\langle\tau+\xi^2\rangle^{b_1}\widehat{v}(\tau,\xi)$ and using
duality, we deduce that Estimate \eqref{prop-g1.1} is equivalent to
\begin{equation} \label{prop-g1.3}
I \lesssim \|f\|_{L^2_{\tau, \xi}}\|g\|_{L^2_{\tau,
\xi}}\|h\|_{L^2_{\tau, \xi}},
\end{equation}
where
\begin{equation} \label{prop-g1.4}
I=\int_{\mathbb R^4} \frac{h(\tau, \xi)\langle
\xi_1\rangle^{1/2}f(\tau_1, \xi_1)g(\tau_2, \xi_2)}{\langle
\sigma\rangle^{c_1}\langle \sigma_1\rangle^{b}\langle
\sigma_2\rangle^{b_1}}d\xi d\xi_1d\tau d\tau_1,
\end{equation}
with $\xi_2=\xi-\xi_1$, $\tau_2=\tau-\tau_1$, $\sigma=\tau+\xi^2$,
$\sigma_1=\tau_1+\gamma|\xi_1|\xi_1$ and $\sigma_2=\tau_2+\xi_2^2$.
The algebraic relation associated to \eqref{prop-g1.4} is given by
\begin{equation} \label{prop-g1.5}
-\sigma+\sigma_1+\sigma_2=-\xi^2+\gamma|\xi_1|\xi_1+\xi_2^2.
\end{equation}

We split the integration domain $\mathbb R^4$ in the following
regions
\begin{eqnarray*}
\mathcal{A} &=& \{(\tau,\tau_1,\xi,\xi_1) \in \mathbb{R}^4 \ : \
|\xi_1|\le 1 \},   \\
\mathcal{B} &=& \{(\tau,\tau_1,\xi,\xi_1) \in \mathbb{R}^4\ :
|\xi_1|> 1\ \mbox{and} \
|\sigma_1|=\max(|\sigma|,|\sigma_1|,|\sigma_2|)\}, \\
\mathcal{C}&=& \{(\tau,\tau_1,\xi,\xi_1) \in \mathbb{R}^4\ :
|\xi_1|> 1\ \mbox{and} \
|\sigma|=\max(|\sigma|,|\sigma_1|,|\sigma_2|)\}, \\
\mathcal{D} &=& \{(\tau,\tau_1,\xi,\xi_1) \in \mathbb{R}^4\ :
|\xi_1|> 1\ \mbox{and} \
|\sigma_2|=\max(|\sigma|,|\sigma_1|,|\sigma_2|)\},
\end{eqnarray*}
and denote by $I_{\mathcal{A}}$, $I_{\mathcal{B}}$,
$I_{\mathcal{C}}$ and $I_{\mathcal{D}}$ the restriction of $I$ to
each one of these regions.

\noindent \textit{Estimate for $I_{\mathcal{A}}$.} In this region
$\langle \xi_1 \rangle \le 1$, then we deduce using Plancherel's
identity and H\"older's inequality that {\setlength\arraycolsep{2pt}
\begin{eqnarray*}
I_{\mathcal{A}} &\lesssim&
\int_{\mathbb{R}^2}\left(\frac{h(\tau,\xi)}{\langle\tau+\xi^2\rangle^{c_1}}\right)^{\vee}
\left(\frac{f(\tau,\xi)}{\langle\tau+\gamma|\xi|\xi\rangle^{b}}\right)^{\vee}
\left(\frac{g(\tau,\xi)}{\langle\tau+\xi^2\rangle^{b_1}}\right)^{\vee}dtdx
\nonumber \\
&\lesssim&\|\left(\frac{h(\tau,\xi)}{\langle\tau+\xi^2\rangle^{c_1}}\right)^{\vee}\|_{L^3_{t,x}}
\|\left(\frac{f(\tau,\xi)}{\langle\tau+\gamma|\xi|\xi\rangle^{b}}\right)^{\vee}\|_{L^3_{t,x}}
\|\left(\frac{g(\tau,\xi)}{\langle\tau+\xi^2\rangle^{b_1}}\right)^{\vee}\|_{L^3_{t,x}}.
\end{eqnarray*}}
This implies, together with Estimates \eqref{prop-g2.1} and
\eqref{prop-g2.2}, that
\begin{equation} \label{prop-g1.6}
I_{\mathcal{A}} \lesssim \|f\|_{L^2_{\tau, \xi}}\|g\|_{L^2_{\tau,
\xi}}\|h\|_{L^2_{\tau, \xi}},
\end{equation}
since $b, \ b_1, \ c_1 > \frac14$.

\noindent \textit{Estimate for $I_{\mathcal{B}}$.} Using the
Cauchy-Schwarz inequality two times, we deduce that
\begin{equation} \label{prop-g1.7}
I_{\mathcal{B}} \lesssim \left(\sup_{\xi_1, \sigma_1}\langle
\sigma_1 \rangle^{-2b}\int_{\mathbb R^2} \frac{|\xi_1|}{\langle
\sigma \rangle^{2c_1}\langle \sigma_2 \rangle^{2b_1}}d\xi
d\sigma\right)^{\frac12} \|f\|_{L^2_{\tau,
\xi}}\|g\|_{L^2_{\tau,\xi}}\|h\|_{L^2_{\tau, \xi}}.
\end{equation}
Remembering the algebraic relation \eqref{prop-g1.5}, we have for $
\xi_1, \ \sigma, \ \sigma_1$ fixed that $d\sigma_2=-2\xi_1d\xi$.
Thus we obtain, by change of variables in the inner integral of the
right-hand side of \eqref{prop-g1.7},
\begin{eqnarray*}
\lefteqn{\langle \sigma_1 \rangle^{-2b}\int_{\mathbb R^2}
\frac{|\xi_1|}{\langle \sigma \rangle^{2c_1}\langle \sigma_2
\rangle^{2b_1}}d\xi d\sigma} \\ && \lesssim \langle \sigma_1
\rangle^{-2b} \left(\int_{|\sigma|\le
|\sigma_1|}\frac{d\sigma}{\langle \sigma \rangle^{2c_1}}\right)
\left(\int_{|\sigma_2|\le |\sigma_1|}\frac{d\sigma_2}{\langle
\sigma_2 \rangle^{2b_1}}\right) \lesssim \langle \sigma_1
\rangle^{2(1-(b+b_1+c_1))}\lesssim 1,
\end{eqnarray*}
since $b+b_1+c_1 \ge 1$. Combining this estimate with
\eqref{prop-g1.7}, we have
\begin{equation} \label{prop-g1.8}
I_{\mathcal{B}} \lesssim  \|f\|_{L^2_{\tau,
\xi}}\|g\|_{L^2_{\tau,\xi}}\|h\|_{L^2_{\tau, \xi}}.
\end{equation}

\noindent \textit{Estimate for $I_{\mathcal{C}}$.} By the
Cauchy-Schwarz inequality (applied two times) it is sufficient to
bound
\begin{equation} \label{prop-g1.9}
\langle \sigma \rangle^{-2c_1}\int_{\mathbb R^2}
\frac{|\xi_1|}{\langle \sigma_1 \rangle^{2b}\langle \sigma_2
\rangle^{2b_1}}d\xi_1 d\sigma_1
\end{equation}
independently of $\xi$ and $\sigma$ to obtain that
\begin{equation} \label{prop-g1.10}
I_{\mathcal{C}} \lesssim  \|f\|_{L^2_{\tau,
\xi}}\|g\|_{L^2_{\tau,\xi}}\|h\|_{L^2_{\tau, \xi}}.
\end{equation}

Now following \cite{bop1}, we first treat the subregion
$|2((1+\gamma\text{sgn}(\xi_1))\xi_1-\xi)|\ge
\frac{1-|\gamma|}{2}|\xi_1|$. When $\xi, \sigma$ and $\sigma_1$ are
fixed, Identity \eqref{prop-g1.7} implies that
$$d\sigma_2=2((1+\gamma\text{sgn}(\xi_1))\xi_1-\xi)d\xi_1.$$
Hence we deduce that
\begin{eqnarray*}
\lefteqn{\langle \sigma \rangle^{-2c_1}\int_{\mathbb R^2}
\frac{|\xi_1|}{\langle \sigma_1 \rangle^{2b}\langle \sigma_2
\rangle^{2b_1}}d\xi_1 d\sigma_1} \\ && \lesssim \langle \sigma
\rangle^{-2c_1} \left(\int_{|\sigma_1|\le
|\sigma|}\frac{d\sigma_1}{\langle \sigma_1 \rangle^{2b}}\right)
\left(\int_{|\sigma_2|\le |\sigma|}\frac{d\sigma_2}{\langle \sigma_2
\rangle^{2b_1}}\right) \lesssim \langle \sigma_1
\rangle^{2(1-(b+b_1+c_1))}\lesssim 1,
\end{eqnarray*}
since $b+b_1+c_1 \ge 1$.

In the subregion of $\mathcal{C}$ where
$|2((1+\gamma\text{sgn}(\xi_1))\xi_1-\xi)|<
\frac{1-|\gamma|}{2}|\xi_1|$, we have from \eqref{prop-g1.7} that
$$|\xi_1|^2 \lesssim |\xi_1^2+\gamma|\xi_1|\xi_1-2\xi\xi_1|\lesssim |\sigma|.$$
Then, we obtain (by applying Lemma \ref{lem-g2}):
{\setlength\arraycolsep{2pt}
\begin{eqnarray*}
\langle \sigma \rangle^{-2c_1}\int_{\mathbb R^2}
\frac{|\xi_1|}{\langle \sigma_1 \rangle^{2b}\langle \sigma_2
\rangle^{2b_1}}d\xi_1 d\sigma_1 &\lesssim& \langle \sigma
\rangle^{\frac12-2c_1}\int_{\mathbb R}\left(\int_{\mathbb R}\langle
\sigma_1 \rangle^{-2b}\langle \sigma_2
\rangle^{-2b_1}d\sigma_1\right)d\xi_1 \\
&\lesssim& \langle \sigma \rangle^{\frac12-2c_1}\int_{\mathbb R}
\langle \sigma+\xi_1^2+\gamma |\xi_1|\xi_1-2\xi\xi_1
\rangle^{1-2(b+b_1)}d\xi_1.
\end{eqnarray*}}
Performing the change of variable
$y=(\theta\xi_1-\theta^{-1}\xi)^2$, where
$\theta=(1+\text{sgn}(\xi_1)\gamma)^{\frac12}$ and noticing that
$|y|\lesssim |\sigma|$ and $dy=2\theta|y|^{\frac12}d\xi_1$ we deduce
that
\begin{eqnarray*}
\langle \sigma \rangle^{-2c_1}\int_{\mathbb R^2}
\frac{|\xi_1|}{\langle \sigma_1 \rangle^{2b}\langle \sigma_2
\rangle^{2b_1}}d\xi_1 d\sigma_1 &\lesssim& \langle \sigma
\rangle^{\frac12-2c_1}\int_{|y|\lesssim |\sigma|}
\frac{dy}{|y|^{\frac12}\langle y-\theta^{-2}\xi^2+\sigma
\rangle^{-(1-2(b+b_1))}}.
\end{eqnarray*}}
Now we use Lemma \ref{lem-g1} to bound the right-hand side integral
by
$$\int_{|y|\lesssim |\sigma|}
|y|^{-\frac12}\langle y \rangle^{1-(2(b+b_1))}dy \lesssim \langle
\sigma \rangle^{[\frac32-2(b+b_1)]_+},$$ where $[\alpha]_+=\alpha$
if $\alpha>0$, $[\alpha]_+=\epsilon$ arbitrarily small if
$\alpha=0$, and $[\alpha]_+=0$ if $\alpha<0$. Therefore
$$\langle \sigma \rangle^{-2c_1}\int_{\mathbb R^2}
\frac{|\xi_1|}{\langle \sigma_1 \rangle^{2b}\langle \sigma_2
\rangle^{2b_1}}d\xi_1 d\sigma_1 \lesssim \langle \sigma
\rangle^{\frac12-2c_1+[\frac32-2(b+b_1)]_+}$$ which is always
bounded using the assumptions on $b$, $b_1$ and $c_1$.

\noindent \textit{Estimate for $I_{\mathcal{D}}$.} By the
Cauchy-Schwarz method it suffices to bound
\begin{equation} \label{prop-g1.11}
\langle \sigma_2 \rangle^{-2b_1}\int_{\mathbb R^2}
\frac{|\xi_1|}{\langle \sigma_1 \rangle^{2b}\langle \sigma
\rangle^{2c_1}}d\xi_1 d\sigma_1
\end{equation}
independently of $\xi_2$ and $\sigma_2$.

We first treat the subregion
$|2((1-\gamma\text{sgn}(\xi_1))\xi_1+\xi_2)|\ge
\frac{1-|\gamma|}{2}|\xi_1|$. When $\xi_2, \ \sigma_2$ and $\sigma$
are fixed, Identity \eqref{prop-g1.7} implies that
$$d\sigma=2(\xi_1+\xi_2-\gamma\text{sgn}(\xi_1)\xi_1)d\xi_1$$
Thus we can estimate \eqref{prop-g1.11} by
\begin{displaymath}
\langle \sigma_2 \rangle^{-2b_1} \left(\int_{|\sigma_1|\le
|\sigma_2|}\frac{d\sigma_1}{\langle \sigma_1 \rangle^{2b}}\right)
\left(\int_{|\sigma|\le |\sigma_2|}\frac{d\sigma}{\langle \sigma
\rangle^{2c_1}}\right) \lesssim \langle \sigma
\rangle^{2(1-(b+b_1+c_1))},
\end{displaymath}
which is bounded since $b+b_1+c_1\ge 1$.

In the subregion $|2((1-\gamma\text{sgn}(\xi_1))\xi_1+\xi_2)|<
\frac{1-|\gamma|}{2}|\xi_1|$, using Identity \eqref{prop-g1.7}, we
deduce that
$$|\xi_1|^2 \lesssim |\xi_1^2-\gamma|\xi_1|\xi_1+2\xi_2\xi_1|\lesssim |\sigma_2|,
$$
where the implicit constant depends on $\gamma$. Then, Lemma
\ref{lem-g2} implies that {\setlength\arraycolsep{2pt}
\begin{eqnarray*}
\langle \sigma_2 \rangle^{-2b_1}\int_{\mathbb R^2}
\frac{|\xi_1|}{\langle \sigma_1 \rangle^{2b}\langle \sigma
\rangle^{2c_1}}d\xi_1 d\sigma_1 &\lesssim& \langle \sigma_2
\rangle^{\frac12-2b_1}\int_{\mathbb R}\left(\int_{\mathbb R}\langle
\sigma_1 \rangle^{-2b}\langle \sigma
\rangle^{-2c_1}d\sigma_1\right)d\xi_1 \\
&\lesssim& \langle \sigma_2 \rangle^{\frac12-2b_1}\int_{\mathbb R}
\langle \sigma_2+(1-\gamma\text{sgn}(\xi_1))\xi_1^2+2\xi_2\xi_1
\rangle^{1-2(b+c_1)}d\xi_1.
\end{eqnarray*}}
We perform the change of variable
$y=(\theta\xi_1+\theta^{-1}\xi_2)^2$ where
$\theta=(1-\gamma\text{sgn}(\xi_1))^{\frac12}$ in the last integral
and we use Lemma \ref{lem-g1} plus the assumptions on $b, \ b_1$ and
$c_1$ to bound \eqref{prop-g1.11} by
\begin{displaymath}
\langle \sigma_2 \rangle^{\frac12-2b_1} \int_{|y|\lesssim
|\sigma_2|}|y|^{-\frac12}\langle y-\theta^{-2}\xi_2^2+ \sigma_2
\rangle^{1-2(b+c_1)}\lesssim \langle \sigma_2
\rangle^{\frac12-2b_1+[\frac32-2(b+c_1)]_+}\lesssim 1.
\end{displaymath}
Therefore, we deduce that
$$I_{\mathcal{D}} \lesssim  \|f\|_{L^2_{\tau,
\xi}}\|g\|_{L^2_{\tau,\xi}}\|h\|_{L^2_{\tau, \xi}},$$ which combined
with \eqref{prop-g1.6}, \eqref{prop-g1.8} and \eqref{prop-g1.10}
implies \eqref{prop-g1.1}.

The proof of Estimate \eqref{prop-g1.2} is actually identical to
that of Estimate \eqref{prop-g1.1}. Indeed, letting $f(\tau,
\xi)=\langle\tau+\xi^2\rangle^{b_1}\widehat{u}(\tau, \xi)$ and
$g(\tau,\xi)=\langle\tau-\xi^2\rangle^{b_1}\widehat{v}(\tau,\xi)$,
we conclude that \eqref{prop-g1.2} is equivalent to
\begin{equation} \label{prop-g1.12}
J \lesssim \|f\|_{L^2_{\tau, \xi}}\|g\|_{L^2_{\tau,
\xi}}\|h\|_{L^2_{\tau, \xi}},
\end{equation}
where
\begin{displaymath}
J=\int_{\mathbb R^4} \frac{|\xi|\langle \xi \rangle^{\frac12}h(\tau,
\xi)f(\tau_1, \xi_1)g(\tau_2, \xi_2)}{\langle
\sigma\rangle^{c}\langle \sigma_1\rangle^{b_1}\langle
\sigma_2\rangle^{b_1}}d\xi d\xi_1d\tau d\tau_1,
\end{displaymath}
with $\xi_2=\xi-\xi_1$, $\tau_2=\tau-\tau_1$,
$\sigma=\tau+\gamma|\xi|\xi$, $\sigma_1=\tau_1+\xi_1^2$ and
$\sigma_2=\tau_2-\xi_2^2$. The algebraic relation associated to this
integral is given by
\begin{equation} \label{prop-g1.5}
-\sigma+\sigma_1+\sigma_2=-\gamma|\xi|\xi+\xi_1^2-\xi_2^2.
\end{equation}
Then, we note that Estimate \eqref{prop-g1.12} is exactly the same
as Estimate \eqref{prop-g1.3}, replacing $c$ by $c_1$ and $b_1$ by
$b$, so we have to ask $2b_1+c \ge 1$ instead of $b+b_1+c_1\ge 1$.
\hfill $\square$

We now slightly modify the bilinear estimates of Proposition
\ref{prop-g1}.
\begin{cor} \label{cor-g1}
Let $\gamma \in \mathbb R$ such that $|\gamma|\neq 1$ and $\gamma
\neq 0$. Then, we have for any $\frac14<b, \ b_1, \ c, \ c_1
<\frac12$ and $s\ge 0$.
\begin{equation} \label{cor-g1.1}
\|uv\|_{X^{s,-c_1}} \lesssim
\|u\|_{Y^{s-\frac12,b}_{\gamma}}\|v\|_{X^{0,b_1}}+
\|u\|_{Y^{-\frac12,b}_{\gamma}}\|v\|_{X^{s,b_1}}, \quad \text{if}
\quad b+b_1+c_1 \ge 1,
\end{equation}
\begin{equation} \label{cor-g1.2}
\|\partial_x(u\bar{v})\|_{Y^{s-\frac12,-c}_{\gamma}} \lesssim
\|u\|_{X^{s,b_1}}\|v\|_{X^{0,b_1}}+\|u\|_{X^{0,b_1}}\|v\|_{X^{s,b_1}},
\quad \text{if} \quad 2b_1+c \ge 1.
\end{equation}
\end{cor}

\noindent \textbf{Proof.} For all $s\ge 0$, we have from the
triangle inequality $\langle \xi \rangle^s \lesssim \langle \xi_1
\rangle^s+\langle \xi-\xi_1 \rangle^s$. Thus we obtain, denoting
$\left(J^s\phi\right)^{\wedge}(\xi) = \langle \xi
\rangle^s\widehat{\phi}(\xi)$ and using \eqref{prop-g1.1}, that
{\setlength\arraycolsep{2pt}
\begin{eqnarray*} \|uv\|_{X^{s,-c_1}}
&\lesssim&
\|J^suv\|_{X^{s,-c_1}}+\|uJ^sv\|_{X^{s,-c_1}} \\
&\lesssim&\|J^su\|_{Y^{-\frac12,b}_{\gamma}}\|v\|_{X^{0,b_1}}+
\|u\|_{Y^{-\frac12,b}_{\gamma}}\|J^sv\|_{X^{0,b_1}} \\
&=&\|u\|_{Y^{s-\frac12,b}_{\gamma}}\|v\|_{X^{0,b_1}}+
\|u\|_{Y^{-\frac12,b}_{\gamma}}\|v\|_{X^{s,b_1}}.
\end{eqnarray*}}
This proves Estimate \eqref{cor-g1.1}. Estimate \eqref{cor-g1.2}
follows using similar arguments with \eqref{prop-g1.2} instead of
\eqref{prop-g1.1}. \hfill $\square$

Finally, we conclude this subsection with the proof of
theorem~\ref{thm-g1}:

\medskip

\noindent \textbf{Proof of Theorem \ref{thm-g1}}

\medskip

\noindent \textit{Case $s=0$.} The system \eqref{eq:lsws} is, at
least formally, equivalent to the integral system
\begin{equation} \label{thm-g1.1}
\left\{ \begin{array}{lll}
         u(t):=F_T^1(u,v) &=& \psi_T U(t)u_0-i\alpha\psi_T\int_0^tU(t-t')u(t')v(t')dt', \\
         v(t):=F_T^2(u) &=&
         \psi_TV_{\gamma}(t)v_0+\beta\psi_T\int_0^tV_{\gamma}(t-t')\partial_x(|u(t')|^2)dt'.
         \end{array}
\right.
\end{equation}
Let $(u_0,v_0) \in L^2(\mathbb R)\times H^{-\frac12}(\mathbb R)$, we
want to use a contraction argument to solve \eqref{thm-g1.1} in a
product of balls
\begin{equation} \label{thm-g1.2}
X^{0,b_1}(a_1)\times Y^{-\frac12,b}_{\gamma}(a_2)=\{(u,v) \in
X^{0,b_1}\times Y^{-\frac12,b}_{\gamma} \ : \ \|u\|_{X^{0,b_1}}\le
a_1, \ \|v\|_{Y^{-\frac12,b}} \le a_2\}.
\end{equation}
Estimates \eqref{lem-g3.1}, \eqref{lem-g3.2}, \eqref{lem-g4.1},
\eqref{lem-g4.3}, \eqref{prop-g1.1} and \eqref{prop-g1.2} imply that
\begin{displaymath}
\|F_T^1(u,v)\|_{X^{0,b_1}} \lesssim
T^{\frac12-b_1}\|u_0\|_{L^2}+T^{1-c_1-b_1}\|u\|_{X^{0,b_1}}\|v\|_{Y^{-\frac12,b}},
\end{displaymath}
and
\begin{displaymath}
\|F_T^2(u)\|_{X^{-\frac12,b_1}} \lesssim
T^{\frac12-b}\|v_0\|_{H^{-\frac12}}+T^{1-c-b}\|u\|_{X^{0,b_1}}^2,
\end{displaymath}
for $\frac14 < b, \ b_1, \ c, \ c_1 < \frac12$ such that $b+c\le 1$,
$b_1+c_1\le 1$, $b+b_1+c_1\ge 1$ and $2b_1+c \ge 1$. In the sequel,
we fix $b=b_1=c=c_1=\frac13$. Therefore we deduce taking $a_1\sim
T^{\frac16}\|u_0\|_{L^2}$ and $a_2\sim
T^{\frac16}\|v_0\|_{H^{-\frac12}}$ that $(F_T^1,F_T^2)$ is a
contraction in $X^{0,\frac13}(a_1)\times
Y^{-\frac12,\frac13}_{\gamma}(a_2)$ if and only if
\begin{equation} \label{thm-g1.3}
T^{\frac12}\|v_0\|_{H^{-\frac12}} \lesssim 1,
\end{equation}
and
\begin{equation} \label{thm-g1.4}
T^{\frac12}\|u_0\|_{L^2}^2 \lesssim \|v_0\|_{H^{-\frac12}}.
\end{equation}
This leads to a solution $(u,v)$ of \eqref{eq:lsws} in
$C([0,T];L^2(\mathbb R))\times C([0,T];H^{-\frac12}(\mathbb R))$
satisfying
\begin{equation} \label{thm-g1.5}
\|u\|_{X^{0,\frac13}} \lesssim T^{\frac16}\|u_0\|_{L^2} \quad
\text{and} \quad \|v\|_{Y^{0,\frac13}_{\gamma}} \lesssim
T^{\frac16}\|v_0\|_{H^{-\frac12}},
\end{equation}
whenever $T$ satisfies \eqref{thm-g1.3} and \eqref{thm-g1.4}. Since
the $L^2$-norm of $u$ is a conserved quantity by the SBO flow, we
can suppose that $\|v_0\|_{H^{-\frac12}} \gg \|u_0\|_{L^2}$,
otherwise we can repeat the above argument and extend the solution
globally in time. Hence Condition \eqref{thm-g1.4} is automatically
satisfied and Condition \eqref{thm-g1.3} implies that the iteration
time $T$ must be $T\sim \|v_0\|_{H^{-\frac12}}^{-2}$. Then we deduce
from \eqref{lem-g4.4}, \eqref{prop-g1.2}, \eqref{thm-g1.1} and
\eqref{thm-g1.5} that there exists a positive constant $C$ such that
\begin{displaymath}
\|v(T)\|_{H^{-\frac12}} \le
\|v_0\|_{H^{-\frac12}}+CT^{\frac12}\|u_0\|_{L^2}^2,
\end{displaymath}
so that we obtain after $m$ iterations of time $T$ where $m \sim
\frac{\|v_0\|_{H^{-\frac12}}}{T^{\frac12}\|u_0\|_{L^2}^2}$ that
\begin{displaymath}
\|v(\Delta T)\|_{H^{-\frac12}}=\|v(mT)\|_{H^{-\frac12}} \le
2\|v_0\|_{H^{-\frac12}}, \quad \text{where} \quad \Delta T \sim
\frac{1}{\|u_0\|_{L^2}^2}.
\end{displaymath}
Since $\Delta T$ only depends on $\|u_0\|_{L^2}$, we can repeat the
above argument and extend the solution $(u,v)$ of \eqref{eq:lsws}
globally in time. Moreover, we deduce that there exists $c>0$ such
that
\begin{equation} \label{thm-g1.6}
\|v(\tilde{T})\|_{H^{-1/2}} \le
e^{c\|u_0\|_{L^2}^2T}\max\{\|u_0\|_{L^2},\|v_0\|_{H^{-1/2}}\}, \
\forall \tilde{T}>0.
\end{equation}

\noindent \textit{Case $s>0$.} Let $(u_0,v_0) \in H^{s}(\mathbb
R)\times H^{s-\frac12}(\mathbb R)$ and $\tilde{T}>0$. This time we
want to solve the integral system \eqref{thm-g1.1} in a space of the
type
\begin{equation} \label{thm-g1.7}
Z^s(a)=\{(u,v) \in X^{s,b_1}\times Y^{s-\frac12,b}_{\gamma} \ : \
\|u\|_{X^{s,b_1}}\le a, \ \|v\|_{Y^{s-\frac12,b}} \le a\}.
\end{equation}
Using Estimates \eqref{lem-g3.1}, \eqref{lem-g3.2},
\eqref{lem-g4.1}, \eqref{lem-g4.3}, \eqref{cor-g1.1} and
\eqref{cor-g1.2} with $b=b_1=c=c_1=\frac13$, we have that for
$0<T\le 1$
\begin{displaymath}
\|F_T^1(u,v)\|_{X^{s,\frac13}} \lesssim
\|u_0\|_{H^s}+T^{\frac13}\left(\|u\|_{X^{0,\frac13}}
\|v\|_{Y^{s-\frac12,\frac13}}
+\|u\|_{X^{s,\frac13}}\|v\|_{Y^{-\frac12,\frac13}}\right),
\end{displaymath}
and
\begin{displaymath}
\|F_T^2(u)\|_{X^{s-\frac12,\frac13}} \lesssim
\|v_0\|_{H^{s-\frac12}}+T^{\frac13}\|u\|_{X^{0,\frac13}}
\|u\|_{X^{s,\frac13}}.
\end{displaymath}
Moreover we can always suppose that $T$ satisfies \eqref{thm-g1.3}
and \eqref{thm-g1.4}, so that Estimate \eqref{thm-g1.5} holds. We
also observe from  the third conservation law in \eqref{eq:law2} and
\textit{a priori} Estimate \eqref{thm-g1.6} that
\begin{equation} \label{thm-g1.8}
\max\{\|u(t)\|_{L^2},\|v(t)\|_{H^{-\frac12}}\}\le
C(\|u_0\|_{L^2},\|v_0\|_{H^{-\frac12}},\tilde{T}), \quad \forall
0<t\le \tilde{T}.
\end{equation}
Therefore we deduce taking
\begin{equation} \label{thm-g1.9}
a\sim\max\{\|u_0\|_{L^2},\|v_0\|_{H^{-1/2}}\} \quad \text{and} \quad
T \sim C(\|u_0\|_{L^2},\|v_0\|_{H^{-\frac12}},\tilde{T})^{-2},
\end{equation} the existence of a unique solution $(u,v)$ of
\eqref{thm-g1.1} in $Z^s(a)$ satisfying the additional regularity $$
(u,v) \in C([0,T];H^s(\mathbb R))\times
C([0,T];H^{s-\frac12}(\mathbb R)).$$ Since the time iteration $T$ in
\eqref{thm-g1.9} only depends on $\|u_0\|_{L^2}$,
$\|v_0\|_{H^{-\frac12}}$ and $\tilde{T}$, we can reapply the above
argument a finite number of times and extend the solution $(u,v)$ on
the time interval $[0,\tilde{T}]$. This completes the proof of
Theorem \ref{thm-g1} if one remembers that $\tilde{T}>0$ is
arbitrary. \hfill $\square$

\setcounter{section}{3} \setcounter{equation}{0} \indent
\begin{center}
{3.2 \ \it Global well-posedness on $\mathbb{T}$}
\end{center}

This subsection contains sharp bilinear estimates for the coupling
terms $uv$ and $\partial_x(|u|^2)$ of the SBO system in the periodic
setting and the global well-posedness result in the energy space
$H^1(\mathbb{T})\times H^{1/2}(\mathbb{T})$ (which is necessary for
our subsequent stability theory).

Let us state our well-posedness result:

\begin{thm}[Local well-posedness in $\mathbb{T}$]\label{theo1}
Let $\gamma \in \Bbb R$ such that $\gamma \neq 0$, $|\gamma| \neq 1$
and $s \ge 1/2$. Then, the SBO system (\ref{eq:lsws}) is locally
well-posed in $H^s(\mathbb{T}) \times H^{s-1/2}(\mathbb{T})$,
$\textit{i.e.}$ for all $(u_0,v_0) \in H^s(\Bbb T)\times
H^{s-1/2}(\Bbb T)$, there exists
$T=T(\|u_0\|_{H^s},\|v_0\|_{H^{s-1/2}})$ and a unique solution of
the Cauchy problem (\ref{eq:lsws}) of the form $(\psi_T u,\psi_T v)$
such that $(u,v) \in X^{s,1/2+}_{per}\times
Y^{s-1/2,1/2+}_{\gamma,per}$. Moreover, $(u,v)$ satisfies the
additional regularity
\begin{equation} \label{theo1.1}
(u,v) \in C([0,T];H^s(\T))\times C([0,T];H^{s-1/2}(\T))
\end{equation}
and the map solution $S: (u_0,v_0) \mapsto (u,v)$ is smooth.
\end{thm}

Using the conservation laws (\ref{eq:law2}) as in \cite{per}, our
local existence result implies
\begin{thm}[Global well-posedness in $\mathbb{T}$]\label{theog}
Let $\alpha, \beta, \gamma \in \Bbb R$ such that $\gamma \neq 0$,
$|\gamma| \neq 1$ and $\frac{\alpha\gamma}{\beta}<0$. Then, the SBO
system (\ref{eq:lsws}) is globally well-posed in $H^s(\mathbb{T})
\times H^{s-1/2}(\mathbb{T})$, when $s \ge 1$.
\end{thm}

The fundamental technical points in the proof of Theorem \ref{theo1}
are the following bilinear estimates. The rest of the proof follows
by standard arguments, as in \cite{CKSTT}.
\begin{prop} \label{prop1}
Let $\gamma \in \Bbb R$ such that $\gamma \neq 0$, $|\gamma| \neq 1$
and $s \ge 1/2$. Then
\begin{equation} \label{prop1.1}
\|uv\|_{X_{per}^{s,-1/2+}} \lesssim
\|u\|_{Y_{\gamma,per}^{s-1/2,1/2}} \|v\|_{X_{per}^{s,1/2}},
\end{equation}
\begin{equation} \label{prop1.2}
\|\partial_x(u\bar{v})\|_{Y^{s-1/2,-1/2+}_{\gamma,per}} \lesssim
\|u\|_{X_{per}^{s,1/2}}\|v\|_{X_{per}^{s,1/2}},
\end{equation}
where the implicit constants depend on $\gamma$.
\end{prop}
These estimates are sharp in the following sense
\begin{prop} \label{prop2}
Let $\gamma\neq 0$, $|\gamma|\neq1$. Then
\begin{itemize}
\item[(i)]The estimate (\ref{prop1.1}) fails for any $s<1/2$.
\item[(ii)]The estimate (\ref{prop1.2}) fails for any $s<1/2$.
\end{itemize}
\end{prop}
\begin{prop} \label{prop3}
Let $\gamma \in \mathbb{R}$ such that $|\gamma|=1$. Then
\begin{itemize}
\item[(i)]The estimate (\ref{prop1.1}) fails for any $s \in \mathbb{R}$.
\item[(ii)]The estimate (\ref{prop1.2}) fails for any $s \in \mathbb{R}$.
\end{itemize}
\end{prop}

The following Bourgain-Strichartz estimates will be used in the
proof of Proposition~\ref{prop1}:
\begin{prop} \label{prop4}
We have
\begin{equation} \label{prop4.1}
\|u\|_{L^4_{t,x}} \lesssim \|u\|_{X^{0,3/8}},
\end{equation}
and \begin{equation} \label{prop4.2}
 \|u\|_{L^4_{t,x}} \lesssim
\|u\|_{Y_{\gamma}^{0,3/8}},
\end{equation}
for $u: \mathbb{R} \times \mathbb{T} \rightarrow \mathbb{C}$ and
$\gamma \in \mathbb{R}$, $\gamma \neq 0$.
\end{prop}

\noindent \textbf{Proof.} The first estimate of (\ref{prop4.1}) was
proved by Bourgain in \cite{Bo} and the second one is a simple
consequence of the first one (see for example \cite{Mo}). \hfill $\square$\\

\noindent \textbf{Proof of Proposition \ref{prop1}.} Fix $s \ge
1/2$. Without loss of generality we can suppose that $0<|\gamma|<1$
in the rest of the proof.

In order to prove estimate (\ref{prop1.1}), it is sufficient to
prove that
\begin{equation} \label{prop1.4}
\|uv\|_{X_{per}^{s,-3/8}}=\|\frac{\langle
n\rangle^s}{\langle\tau+n^2\rangle^{3/8}}(uv)^{\wedge}
(\tau,n)\|_{l^2_nL^2_{\tau}}\lesssim
\|u\|_{Y_{\gamma,per}^{s-1/2,1/2}}\|v\|_{X_{per}^{s,1/2}}.
\end{equation}

Letting $f(\tau, n)=\langle n\rangle^{s-1/2}\langle\tau+\gamma
n|n|\rangle^{1/2}\widehat{u}(\tau, n)$, $g(\tau, n)=\langle
n\rangle^{s}\langle\tau+n^2\rangle^{1/2}\widehat{v}(\tau, n)$ and
using duality, we deduce that Estimate (\ref{prop1.4}) is equivalent
to
\begin{equation} \label{prop1.5}
I \lesssim
\|f\|_{L^2_{\tau}l^2_n}\|g\|_{L^2_{\tau}l^2_n}\|h\|_{L^2_{\tau}l^2_n},
\end{equation}
where {\setlength\arraycolsep{2pt}
\begin{eqnarray}  \label{prop1.6}
I &:=& \sum_{n, n_1\in\mathbb{Z}}\int_{\mathbb{R}^2}\frac{\langle
n\rangle^s h(\tau, n)f(\tau_1, n_1)}
{\langle\tau+n^2\rangle^{3/8}\langle n_1\rangle^{s-1/2}\langle\tau_1+\gamma|n_1|n_1\rangle^{1/2}} \nonumber \\
&& \times \frac{g(\tau-\tau_1, n-n_1)}{\langle n-n_1\rangle^{s}
\langle\tau-\tau_1+(n-n_1)^2\rangle^{1/2}}d\tau d\tau_1.
\end{eqnarray}}
In order to bound the integral in (\ref{prop1.6}), we split the
integration domain $\mathbb{R}^2\times \mathbb{Z}^2$ in the
following regions,
\begin{eqnarray*}
\mathcal{M} &=& \{(\tau,\tau_1,n,n_1) \in
\mathbb{R}^2\times\mathbb{Z}^2 \ : \ n_1=0 \
\mbox{or} \ |n| \le c(\gamma)^{-1}|n-n_1|\},   \\
\mathcal{N} &=& \{(\tau,\tau_1,n,n_1) \in
\mathbb{R}^2\times\mathbb{Z}^2 \ : \ n_1 \neq 0 \ \mbox{and} \
|n-n_1| \le c(\gamma)|n|\},
\end{eqnarray*}
where $c(\gamma)$ is a positive constant depending on $\gamma$ to be
fixed later. We also denote by $I_{\mathcal{M}}$ and
$I_{\mathcal{N}}$ the integral $I$ restricted to the regions
$\mathcal{M}$ and $\mathcal{N}$, respectively.

\noindent \textit{Estimate on the region $\mathcal{M}$.} We observe
that, since $s\geq 1/2$, it holds $\frac{\langle n\rangle^s}{\langle
n_1\rangle^{s-1/2}\langle n-n_1\rangle^s} \lesssim 1$ (where the
implicit constant depends on $\gamma$) in the region $\mathcal{M}$.
Thus, we deduce that, using the Plancherel identity and the
$L^2_{t,x}L^4_{t,x}L^4_{t,x}$-H\"older inequality,
{\setlength\arraycolsep{2pt}
\begin{eqnarray}
I_{\mathcal{M}} & \lesssim &  \sum_{n,
n_1\in\mathbb{Z}}\int_{\mathbb{R}^2}\frac{h(\tau, n)f(\tau_1,
n_1)g(\tau-\tau_1, n-n_1)d\tau d\tau_1}
{\langle\tau+n^2\rangle^{3/8}\langle\tau_1+\gamma|n_1|n_1\rangle^{1/2}\langle\tau-\tau_1+(n-n_1)^2\rangle^{1/2}}
\nonumber \\
&\lesssim& \int_{\mathbb{R} \times
\mathbb{T}}\left(\frac{h(\tau,n)}{\langle\tau+n^2\rangle^{3/8}}\right)^{\vee}
\left(\frac{f(\tau,n)}{\langle\tau+\gamma|n|n\rangle^{1/2}}\right)^{\vee}
\left(\frac{g(\tau,n)}{\langle\tau+n^2\rangle^{1/2}}\right)^{\vee}dtdx
\nonumber \\
&\lesssim&\|\left(\frac{h(\tau,n)}{\langle\tau+n^2\rangle^{3/8}}\right)^{\vee}\|_{L^2_{t,x}}
\|\left(\frac{f(\tau,n)}{\langle\tau+\gamma|n|n\rangle^{1/2}}\right)^{\vee}\|_{L^4_{t,x}}
\|\left(\frac{g(\tau,n)}{\langle\tau+n^2\rangle^{1/2}}\right)^{\vee}\|_{L^4_{t,x}}.
\label{prop1.6b}
\end{eqnarray}}
This implies, together with Estimates (\ref{prop4.1}) and
(\ref{prop4.2}), that
\begin{equation} \label{prop1.6b2}
I_{\mathcal{M}} \lesssim
\|f\|_{L^2_{\tau}l^2_n}\|g\|_{L^2_{\tau}l^2_n}\|h\|_{L^2_{\tau}l^2_n}.
\end{equation}

\noindent \textit{Estimate on the region $\mathcal{N}$.} The
dispersive smoothing effect associated to the SBO system
(\ref{eq:lsws}) can be translated by the following algebraic
relation
\begin{equation} \label{prop1.7}
-(\tau+n^2)+(\tau_1+\gamma n_1|n_1|)+(\tau-\tau_1+(n-n_1)^2)=
Q_{\gamma}(n,n_1),
\end{equation}
where
\begin{equation} \label{prop1.8}
Q_{\gamma}(n,n_1)=(n-n_1)^2+\gamma |n_1|n_1-n^2.
\end{equation}
We have in the region $\mathcal{N}$, $|n_1| \le (1+c(\gamma))|n|$,
so that
\begin{displaymath}
|Q_{\gamma}(n,n_1)| \ge
\left(1-|\gamma|(1+c(\gamma))^2-c(\gamma)^2\right)(1+c(\gamma))^{-2}|n_1|^2.
\end{displaymath}
Now, we choose $c(\gamma)$ positive, small enough such that
$$\left(1-|\gamma|(1+c(\gamma))^2-c(\gamma)^2\right)=\frac{1-|\gamma|}{2},$$
which is possible since $|\gamma|<1$. Therefore, we divide the
region $\mathcal{N}$ in three parts accordingly to which term of the
left-hand side of (\ref{prop1.7}) is dominant:
\begin{eqnarray*}
\mathcal{N}_1 &=& \{(\tau,\tau_1,n,n_1) \in \mathcal{N} \ : \
|\tau+n^2|\ge|\tau_1+\gamma
|n_1|n_1|,|\tau-\tau_1+(n-n_1)^2|\},   \\
\mathcal{N}_2 &=& \{(\tau,\tau_1,n,n_1) \in \mathcal{N} \ : \
|\tau_1+\gamma |n_1|n_1| \ge
|\tau+n^2|,|\tau-\tau_1+(n-n_1)^2| \}, \\
\mathcal{N}_3 &=& \{(\tau,\tau_1,n,n_1) \in \mathcal{N} \ : \
|\tau-\tau_1+(n-n_1)^2|\ge |\tau+n^2|, |\tau_1+\gamma |n_1|n_1| \}.
\end{eqnarray*}
We denote by $I_{\mathcal{N}_1}$, $I_{\mathcal{N}_2}$ and
$I_{\mathcal{N}_3}$ the restriction of the integral $I$ to the
regions $\mathcal{N}_1$, $\mathcal{N}_2$ and $\mathcal{N}_3$,
respectively.

In the region $\mathcal{N}_1$, we have $\frac{\langle
n\rangle^s}{\langle n_1\rangle^{s-1/2}\langle n-n_1\rangle^s} \times
\frac{1}{\langle\tau+n^2\rangle^{3/8}} \lesssim 1$, so that we can
conclude
\begin{equation} \label{prop1.9}
I_{\mathcal{N}_1} \lesssim
\|f\|_{L^2_{\tau}l^2_n}\|g\|_{L^2_{\tau}l^2_n}\|h\|_{L^2_{\tau}l^2_n},
\end{equation}
exactly as in (\ref{prop1.6b}). We note that $\frac{\langle
n\rangle^s}{\langle n_1\rangle^{s-1/2}\langle n-n_1\rangle^s} \times
\frac{1}{\langle\tau_1+|n_1|n_1\rangle^{1/2}} \lesssim 1$ in the
region $\mathcal{N}_2$. Then, using the
$L^4_{t,x}L^2_{t,x}L^4_{t,x}$-H\"older inequality, that
$$I_{\mathcal{N}_2} \lesssim \|\left(\frac{h(\tau,n)}
{\langle\tau+n^2\rangle^{3/8}}\right)^{\vee}\|_{L^4_{t,x}}
\|f\|_{L^2_{\tau}l^2_n}
\|\left(\frac{g(\tau,n)}{\langle\tau+n^2\rangle^{1/2}}\right)^{\vee}\|_{L^4_{t,x}}.$$
Combining this with Estimate (\ref{prop4.1}), we obtain that
\begin{equation} \label{prop1.10}
I_{\mathcal{N}_2} \lesssim
\|f\|_{L^2_{\tau}l^2_n}\|g\|_{L^2_{\tau}l^2_n}\|h\|_{L^2_{\tau}l^2_n}.
\end{equation}
Similarly, $\frac{\langle n\rangle^s}{\langle
n_1\rangle^{s-1/2}\langle n-n_1\rangle^s} \times
\frac{1}{\langle\tau-\tau_1+(n-n_1)^2\rangle^{1/2}} \lesssim 1$ in
$\mathcal{N}_3$ so that {\setlength\arraycolsep{2pt}
\begin{eqnarray}
I_{\mathcal{N}_3} & \lesssim &
\|\left(\frac{h(\tau,n)}{\langle\tau+n^2\rangle^{3/8}}\right)^{\vee}\|_{L^4_{t,x}}
\|\left(\frac{f(\tau,n)}{\<\tau+\gamma|n|n\>^{1/2}}\right)^{\vee}\|_{L^4_{t,x}}
\|g\|_{L^2_{\tau}l^2_n} \nonumber \\
&\lesssim&
\|f\|_{L^2_{\tau}l^2_n}\|g\|_{L^2_{\tau}l^2_n}\|h\|_{L^2_{\tau}l^2_n}.
\label{prop1.11}
\end{eqnarray}}

Then, we gather (\ref{prop1.6b2}), (\ref{prop1.9}), (\ref{prop1.10})
and (\ref{prop1.11}) to deduce (\ref{prop1.5}), which concludes the
proof of the estimate (\ref{prop1.1}).

Next, in order to prove Estimate (\ref{prop1.2}), we argue as above
so that it is sufficient to prove
\begin{equation} \label{prop1.12}
\|\partial_x(u\bar{v})\|_{Y_{\gamma, per}^{s-1/2,-3/8}}\lesssim
\|u\|_{X^{s,1/2}_{per}}\|v\|_{X^{s,1/2}_{per}} ,
\end{equation}
which is equivalent by duality and after performing the change of
variable $f(\tau,n)=\<n\>^s\<\tau+n^2\>^{1/2}\widehat{u}(\tau,n)$
and $g(\tau,n)=\<n\>^s\<\tau-n^2\>^{1/2}\widehat{\bar{v}}(\tau,n)$
to
\begin{equation} \label{prop1.13}
J \lesssim
\|f\|_{L^2_{\tau}l^2_n}\|g\|_{L^2_{\tau}l^2_n}\|h\|_{L^2_{\tau}l^2_n},
\end{equation}
where {\setlength\arraycolsep{2pt}
\begin{eqnarray}  \label{prop1.14}
J &:=& \sum_{n, n_1\in\mathbb{Z}}\int_{\R^2}\frac{|n|\langle
n\rangle^{s-1/2}h(\tau, n)f(\tau_1, n_1)}
{\<\tau+\gamma|n|n\>^{3/8}\<n_1\>^s\<\tau_1+n_1^2\>^{1/2}} \nonumber \\
&& \times \frac{g(\tau-\tau_1, n-n_1)}{\<n-n_1\>^{s}
\<\tau-\tau_1-(n-n_1)^2\>^{1/2}}d\tau d\tau_1.
\end{eqnarray}}
The algebraic relation associated to (\ref{prop1.14}) is given by
$$-(\tau+\gamma |n|n)+(\tau_1+n_1^2)+(\tau-\tau_1-(n-n_1)^2)=
\tilde{Q}_{\gamma}(n,n_1),$$ where
$$\tilde{Q}_{\gamma}(n,n_1)=-(n-n_1)^2-\gamma|n|n+n_1^2.$$
Therefore we can prove Estimate (\ref{prop1.13}) using exactly
the same arguments as for Estimate (\ref{prop1.5}).\hfill $\square$ \\

\noindent \textbf{Remark 3.1.} Observe that we obtain our bilinear
estimates in the spaces $X^{s,1/2+}_{per}$ and
$Y^{s-1/2,1/2+}_{\gamma,per}$ which control the $L^{\infty}_tH^s_x$
and $L^{\infty}_tH^{s-1/2}_x$ norms respectively. Therefore, we do
not need to use other norms as in the case of the periodic KdV
equation \cite{CKSTT}. \\

\noindent\textbf{Remark 3.2.} Observe that the proof of Proposition
\ref{prop1} actually shows that the following bilinear estimates
hold:
\begin{equation*}
\|uv\|_{X_{per}^{s,-3/8}} \lesssim
\|u\|_{X_{per}^{s,3/8}}\|v\|_{Y_{\gamma,per}^{s-1/2,1/2}} +
\|u\|_{X_{per}^{s,1/2}}\|v\|_{Y_{\gamma,per}^{s-1/2,3/8}},
\end{equation*}
\begin{equation*}
\|\partial_x(u\bar{w})\|_{Y^{s-1/2,-3/8}_{\gamma,per}} \lesssim
\|u\|_{X_{per}^{s,3/8}}\|w\|_{X_{per}^{s,1/2}} +
\|u\|_{X_{per}^{s,1/2}}\|w\|_{X_{per}^{s,3/8}},
\end{equation*}
While we are not attempting to use this refined version of
proposition \ref{prop1} in this paper, we plan to apply these
estimates combined with the I-method of Colliander, Keel,
Stafillani, Takaoka and Tao to get global well-posedness results for
the periodic SBO system below the energy space. Indeed, this issue
will be addressed in a forthcoming paper.

\bigskip

In the proof of Proposition \ref{prop2}, we will use the following
lemma which is a direct consequence of the Dirichlet theorem.

\begin{lem} \label{lemm1}
Let $\gamma \in \R$ such that $\gamma \neq 0$ and $|\gamma| < 1$ and
$Q_{\gamma}$ defined as in (\ref{prop1.8}). Then, there exists a
sequence of positive integers $\{N_j\}_{j \in \N}$ such that
\begin{equation} \label{lemm1.1}
N_j \rightarrow \infty \quad \mbox{and} \quad
|Q_{\gamma}(N_j,N_j^0)|\le 1,
\end{equation}
where $N^0_j=[\frac{2N_j}{1+\gamma}]$ and $[x]$ denotes the closest
integer to $x$.
\end{lem}

\begin{thm} [Dirichlet] \label{theo2}
Let $\alpha \in \R \setminus \Q$. Then, the inequality
\begin{equation} \label{theo1.1}
0<\left|\alpha-\frac{p}{q}\right|<\frac{1}{q^2}
\end{equation}
has infinitely many rational solutions $\frac{p}{q}$.
\end{thm}

\noindent \textbf{Proof of Lemma \ref{lemm1}.} Fix $\gamma \in \R$
such that $|\gamma|<1$. Let $N$ a positive integer, $N \ge 2$,
$\alpha=\frac{2}{1+\gamma}$ and $N^0=[\alpha N]$. Then, from the
definition in (\ref{prop1.8}), we deduce that
\begin{equation} \label{lemm1.2}
|Q_{\gamma}(N,N^0)|\le 1 \quad \iff \quad \left|\alpha-\frac{[\alpha
N]}{N}\right|\le \frac{1}{N^2}.
\end{equation}
When $\alpha \in \Q$, $\alpha=\frac{p}{q}$, it is clear that we can
find an infinity of positive integer $N$ satisfying the right-hand
side of (\ref{lemm1.2}) choosing $N_j=jq$, $j \in \N$. When $\alpha
\in \R \setminus \Q$, this is guaranteed by the Dirichlet theorem.
\hfill $\square$ \\

\noindent \textbf{Proof of Proposition \ref{prop2}.} We will only
show that the estimate (\ref{prop1.1}) fails, since a counterexample
for the estimate (\ref{prop1.2}) can be constructed in a similar
way. First observe that, letting $f(\tau,
n)=\<n\>^{s-1/2}\<\tau+\gamma n|n|\>^{1/2}\widehat{u}(\tau, n)$ and
$g(\tau, n)=\<n\>^{s}\<\tau+n^2\>^{1/2}\widehat{v}(\tau, n)$, the
estimate (\ref{prop1.1}) is equivalent to
\begin{equation} \label{prop2.1}
\|B_{\gamma}(f,g;s)\|_{L^2_{\tau}l^2_n} \lesssim
\|f\|_{L^2_{\tau}l^2_n}\|g\|_{L^2_{\tau}l^2_n}, \ \forall \ f, \ g
\in L^2_{\tau}l^2_{n},
\end{equation}
where {\setlength\arraycolsep{2pt} {\begin{eqnarray}
B_{\gamma}(f,g;s)(\tau, n)&:=& \frac{\<n\>^s}{\<\tau+n^2\>^{1/2}}
\sum_{n_1\in\Z}\int_{\R}\frac{f(\tau_1, n_1)}
{\<n_1\>^{s-1/2}\<\tau_1+\gamma|n_1|n_1\>^{1/2}} \nonumber
\\
&& \times \frac{g(\tau-\tau_1,
n-n_1)}{\<n-n_1\>^{s}\<\tau-\tau_1+(n-n_1)^2\>^{1/2}}d\tau_1,
\label{prop2.2}
\end{eqnarray}}
for all $s$ and $\gamma \in \R$.

Fix $s<1/2$ and $\gamma$ such that $|\gamma| \neq 1$; without loss
of generality, we can suppose that $|\gamma|<1$. Consider the
sequence of integer $\{N_j\}$ obtained in Lemma \ref{lemm1}, which
we can always suppose to verify $N_j \gg 1$, and define
\begin{equation} \label{prop2.3}
f_j(\tau,n)=a_n\chi_{1/2}(\tau+\gamma|n|n) \quad \mbox{with} \quad
a_n= \left\{  \begin{array}[pos]{ll}
               1, & n=N^0_j,\\
               0, & \mbox{elsewhere},
              \end{array}   \right.
\end{equation}
and
\begin{equation} \label{prop2.4}
g_j(\tau,n)=b_n\chi_{1/2}(\tau+n^2) \quad \mbox{with} \quad b_n=
\left\{  \begin{array}[pos]{ll}
               1, & n=N_j-N_j^0,\\
               0, & \mbox{elsewhere},
              \end{array}   \right.
\end{equation}
where $\chi_{r}$ is the characteristic function of the interval
$[-r,r]$. Hence,
\begin{equation} \label{prop2.5}
\|f_{j}\|_{L^2_{\tau}l^2_n} \sim  \|g_{j}\|_{L^2_{\tau}l^2_n} \sim
1,
\end{equation}
\begin{displaymath}
a_{n_1}b_{n-n_1} \neq 0 \quad \mbox{if and only if} \quad n_1=N_j^0
\ \mbox{and} \ n=N_j.
\end{displaymath}
Using (\ref{prop1.7}), we deduce that
\begin{displaymath}
\int_{\R}\chi_{1/2}(\tau_1+\gamma|N_j^0|N_j^0)\chi_{1/2}(\tau-\tau_1+(N_j-N_j^0)^2)d\tau_1
\sim \chi_{1}(\tau+N^2+Q_{\gamma}(N_j,N_j^0)).
\end{displaymath}
Therefore, we have from the definition in (\ref{prop2.2})
\begin{equation} \label{prop2.6}
B_{\gamma}(f_j,g_j;s)(\tau,N_j) \gtrsim
\frac{N_j^s\chi_{1}(\tau+N_j^2+Q_{\gamma}(N_j,N_j^0))}{\<\tau+N_j^2\>^{1/2}N_j^{s-1/2}N_j^s},
\end{equation}
where the implicit constant depends on $\gamma$. Thus, we deduce
using (\ref{lemm1.1}) that
\begin{equation} \label{prop2.7}
\|B_{\gamma}(f_j,g_j;s)\|_{L^2_{\tau}l^2_n} \gtrsim N_j^{1/2-s},
\quad \forall \ j \in \N
\end{equation}
which combined with (\ref{lemm1.1}) and (\ref{prop2.5}) contradicts
(\ref{prop2.1}), since $s<1/2$. \hfill $\square$ \\

\noindent \textbf{Proof of Proposition \ref{prop3}.} Let $s \in \R$,
we fix $\gamma = 1$. As in the proof of Proposition \ref{prop2}, we
will only show that the estimate (\ref{prop1.1}) fails, since a
counterexample for the estimate (\ref{prop1.2}) can be constructed
in a similar way. In this case, (\ref{prop1.1}) is equivalent to
\begin{equation} \label{prop3.1}
\|B_{1}(f,g;s)\|_{L^2_{\tau}l^2_n} \lesssim
\|f\|_{L^2_{\tau}l^2_n}\|g\|_{L^2_{\tau}l^2_n}, \ \forall \ f, \ g
\in L^2_{\tau}l^2_{n},
\end{equation}
where {\setlength\arraycolsep{2pt} {\begin{eqnarray}
B_{1}(f,g;s)(\tau, n)&:=& \frac{\<n\>^s}{\<\tau+n^2\>^{1/2}}
\sum_{n_1\in\Z}\int_{\R}\frac{f(\tau_1, n_1)}
{\<n_1\>^{s-1/2}\<\tau_1+|n_1|n_1\>^{1/2}} \nonumber
\\
&& \times \frac{g(\tau-\tau_1,
n-n_1)}{\<n-n_1\>^{s}\<\tau-\tau_1+(n-n_1)^2\>^{1/2}}d\tau_1.
\label{prop3.2}
\end{eqnarray}}
Fix a positive integer $N$, such that $N \gg 1$, and define
\begin{equation} \label{prop3.3}
f_N(\tau,n)=a_n\chi_{1/2}(\tau+|n|n) \quad \mbox{with} \quad a_n=
\left\{  \begin{array}[pos]{ll}
               1, & n=N,\\
               0, & \mbox{elsewhere},
              \end{array}   \right.
\end{equation}
and
\begin{equation} \label{prop3.4}
g_N(\tau,n)=b_n\chi_{1/2}(\tau+n^2) \quad \mbox{with} \quad b_n=
\left\{  \begin{array}[pos]{ll}
               1, & n=0,\\
               0, & \mbox{elsewhere},
              \end{array}   \right.
\end{equation}
where $\chi_{r}$ is the characteristic function of the interval
$[-r,r]$. Hence,
\begin{equation} \label{prop3.5}
\|f_{N}\|_{L^2_{\tau}l^2_n} \sim  \|g_{N}\|_{L^2_{\tau}l^2_n} \sim
1,
\end{equation}
\begin{displaymath}
a_{n_1}b_{n-n_1} \neq 0 \quad \mbox{if and only if} \quad n_1=N \
\mbox{and} \ n=N,
\end{displaymath}
and
\begin{displaymath}
\int_{\R}\chi_{1/2}(\tau_1+N^2)\chi_{1/2}(\tau-\tau_1)d\tau_1 \sim
\chi_{1}(\tau+N^2).
\end{displaymath}
Therefore, we deduce from (\ref{prop3.2}) that
\begin{equation} \label{prop3.6}
\|B_{1}(f_N,g_N;s)\|_{L^2_{\tau}l^2_n} \gtrsim N^{1/2}, \quad
\forall \ N \gg 1,
\end{equation}
which combined with (\ref{prop3.5}) contradicts (\ref{prop3.1}). The
case $\gamma=-1$ is similar. \hfill $\square$

\begin{center}
\item\section{Existence of Periodic Traveling-Wave Solutions}
\end{center}

The goal of this section is to show the existence of a smooth branch
of periodic traveling-wave solutions for (\ref{eq:eds}). Initially
we show a novel smooth branch of dnoidal waves solutions for
(\ref{eq:eds}) in the case $\gamma=0$. After that, by using the
implicit function theorem, we construct (in the case $\gamma\neq 0$)
a smooth curve of periodic solutions bifurcating from these dnoidal
waves.

\setcounter{section}{4} \setcounter{equation}{0}
\begin{center}
{4.1 \ \it Dnoidal Waves Solutions}
\end{center}

We start by finding solutions for the case $\gamma=0$ and $\sigma>0$
in (\ref{eq:eds}).  Henceforth, without loss of generality, we will
assume that $\alpha=1$ and $\beta=\frac12$. Hence, we need to solve
the system
\begin{equation}\label{eq:2.1}
\left \{
    \begin{array}{l}
\phi_0''-\sigma\phi_0=\psi_0\phi_0\\
\psi_0=-\frac{1}{2c}\phi_0^2.
\end{array}
\right.
\end{equation}
Then, by replacing the second equation of (\ref{eq:2.1}) into the
first one, we obtain that $\phi_0$ satisfies
\begin{equation}\label{eq:2.2}
\phi_0''-\sigma\phi_0+\frac{1}{2c}\phi_0^3=0.
\end{equation}
Equation (\ref{eq:2.2}) can be solved in a similar fashion to the
method used by Angulo \cite{a:a} (in the context of periodic
traveling-wave solutions for the nonlinear Schr\"odinger equation
(\ref{sch})). For the sake of completeness, we provide here a sketch
of the proof of this fact. Indeed, from (\ref{eq:2.2}), $\phi_0$
must satisfy the first-order equation
\begin{equation}\label{eq:2.3}
[\phi_0']^2=\frac{1}{4c}[-\phi_0^4+4c\sigma\phi_0^2+4cB_{\phi_0}]\nonumber
=\frac{1}{4c}(\eta_1^2-\phi_0^2)(\phi_0^2-\eta_2^2),
\end{equation}
where $B_{\phi_0}$ is an integration constant and $-\eta_1,
\eta_1, -\eta_2, \eta_2$ are the zeros of the polynomial
$F(t)=-t^4+ 4c\sigma t^2+4cB_{\phi_0}$. Moreover,
\begin{equation}\label{eq:2.4}
\left \{
    \begin{array}{l}
4c\sigma=\eta_1^2+\eta^2_2\\
4cB_{\phi_0}=-\eta_1^2\eta_2^2.
\end{array}
\right.
\end{equation}
We suppose, without loss of generality, that $\eta_1>\eta_2>0$. Then
$\eta_2\leq \phi_0\leq\eta_1$ and so $\phi_0$ will be a positive
solution. Note that $-\phi_0$ is also a solution of (\ref{eq:2.2}).
Next, define $\zeta=\phi_0/\eta_1$ and
$k^2=(\eta_1^2-\eta_2^2)/\eta_1^2$. It follows from (\ref{eq:2.3})
that
\begin{eqnarray*}
[\zeta']^2&=\frac{\eta_1^2}{4c}(1-\zeta^2)(\zeta^2+k^2-1).
\end{eqnarray*}
Let us now define a new function $\chi$ through
$\zeta^2=1-k^2\sin^2\chi$. So we get that $
4c(\chi')^2=\eta_1^2(1-k^2\sin^2\chi)$. Then for
$l=\frac{\eta_1}{2\sqrt c}$, and  assuming that $\zeta(0)=1$, we
have
$$
\int_0^{\chi(\xi)}\frac{dt}{\sqrt{1-k^2sin^2t}}\;=l\;\xi.
$$
Then from the definition of the {\it Jacobian elliptic function}
$sn(u;k)$, we have that $\sin \chi=sn(l\xi;k)$ and hence $
\zeta(\xi)=\sqrt{1-k^2sn^2(l\xi;k)}\equiv dn(l\xi;k)$. Returning to
the variable $\phi_0$, we obtain the novel {\bf{dnoidal waves}}
solutions associated to equation (\ref{eq:2.1}),
\begin{equation}\label{eq:2.5}
\left \{
    \begin{array}{l}
\phi_0(\xi)\equiv
\phi_0(\xi;\eta_1,\eta_2)=\eta_1\;dn\Big(\frac{\eta_1}{2\sqrt{c}}\;\xi;k\Big)\\
\psi_0(\xi)\equiv
\psi_0(\xi;\eta_1,\eta_2)=-\frac{\eta_1^2}{2c}\;dn^2\Big(\frac{\eta_1}{2\sqrt{c}}\;\xi;k\Big),
\end{array}
\right.
\end{equation}
where
\begin{equation}\label{eq:2.6}
0<\eta_2<\eta_1,\qquad
k^2=\frac{\eta_1^2-\eta_2^2}{\eta_1^2},\qquad
\eta_1^2+\eta_2^2=4c\sigma.
\end{equation}
Next, since $dn$ has fundamental period $2K(k)$, it follows that
$\phi_0$ in (\ref{eq:2.5}) has fundamental wavelength (i.e., period)
$T_{\phi_0}$ given by
$$
T_{\phi_0}\equiv \frac{4\sqrt{c}}{\eta_1}\;K(k).
$$
Given $c>0$, $\sigma>0$, it follows from (\ref{eq:2.6}) that $0<
\eta_2<\sqrt {2c\sigma} <\eta_1<2\sqrt {c\sigma}$. Moreover we can
write
\begin{equation}\label{eq:2.8}
T_{\phi_0}(\eta_2)=\frac{4\sqrt{c}}{\sqrt{4c\sigma-\eta_2^2}}\;K(
k(\eta_2))\qquad\text{with}\qquad
k^2(\eta_2)=\frac{4c\sigma-2\eta_2^2}{4c\sigma-\eta_2^2}.
\end{equation}
Then, using these formulas and the properties of the function $K$,
we see that $T_{\phi_0}\in (\sqrt{\frac{2}{\sigma}}\;\pi,+\infty)$
for $\eta_2\in (0, \sqrt{2c\sigma})$. Moreover, we will see in
Theorem \ref{thm2.1} below that $\eta_2\mapsto T_{\phi_0}(\eta_2)$
is a strictly decreasing mapping and so we obtain the basic
inequality
\begin{equation}\label{eq:2.9}
T_{\phi_0}>\sqrt{\frac{2}{\sigma}}\;\pi.
\end{equation}

Two relevant solutions of (\ref{eq:2.1}) are hidden in
(\ref{eq:2.5}). Namely, the constant and solitary wave solutions.
Indeed, when $\eta_2\to\sqrt{2c\sigma}$, i.e. $\eta_2\to \eta_1$, it
follows that $k\to
 0^+$. Then since $d(u;0^+)\to 1$ we obtain the constant solutions
\begin{equation}\label{eq:2.10}
\phi_0(\xi)=\sqrt{2c\sigma}\qquad\text{and}\qquad
\psi_0(\xi)=-\sigma.
\end{equation}
Next, for $\eta_2\to 0$ we have $\eta_1\to 4c\sigma^{-}$ and so
$k\to 1^{-}$. Then since $dn(u;1^{-})\to sech(u)$ we obtain the
classical solitary wave solutions
\begin{equation}\label{eq:2.11}
\phi_{0,s}(\xi)=2\sqrt{c\sigma}sech(\sqrt{\sigma} \xi)
\qquad\text{and}\qquad \psi_{0,s}(\xi)=-2\sigma
sech^2(\sqrt{\sigma} \xi).
\end{equation}

Our next theorem is the main result of this subsection and it proves
that for a fixed period $L>0$ there exists a smooth branch of
dnoidal waves solutions with the same period $L$ to the system
(\ref{eq:2.1}) (or equivalently to equation (\ref{eq:2.2})). The
construction of a family  of dnoidal waves with a fixed period $L$
is an immediate  consequence of the analysis made above. Indeed, let
$L>0$ be a fixed number. Choose $c>0$ and $w\in \Bbb R$ real fixed
numbers such that $\sigma\equiv\omega-c^2/4>2\pi^2/{L^2}$. Since the
function $\eta_2\in (0,\sqrt{2c\sigma})\to T_{\phi_0}(\eta_2)$ is
strictly decreasing (see Theorem \ref{thm2.1} below), there is a
unique $\eta_2=\eta_2(\sigma)\in (0,\sqrt{2c\sigma})$ such that
$\phi_0(\cdot;\eta_1(\sigma),\eta_2(\sigma))$ has fundamental period
$T_{\phi_0}(\eta_2(\sigma))=L$. We claim that the choice of
$\eta_2(\sigma)$ depends smoothly of $\sigma$:

\begin{thm}  Let $L$ and $c$ be arbitrarily fixed positive
numbers. Let $\sigma_0>2\pi^2/L^2$ and $\eta_{2,0}=\eta
_2(\sigma_0)$ be the unique number in the interval $
(0,\sqrt{2c\sigma})$ such that $T_{\phi_0}(\eta_{2,0})=L$. Then,

$(1)$\; there are intervals $I(\sigma_0)$ and $B(\eta_{2,0})$ around
of $\sigma_0$ and $\eta_2(\sigma_0)$ respectively, and an unique
smooth function $\Lambda:I(\sigma_0)\to B(\eta_{2,0})$, such  that
$\Lambda(\sigma_0)= \eta_{2,0}$ and
\begin{equation}\label{eq:2.12}
\frac{4\sqrt c}{\sqrt{4c\sigma-\eta_2^2}}\;K(k(\sigma))=L,
\end{equation}
where $\sigma\in I(\sigma_0)$, $\eta_2=\Lambda(\sigma)$, and
\begin{equation}\label{eq:2.13}
k^2\equiv
k^2(\sigma)=\frac{4c\sigma-2\eta_2^2}{4c\sigma-\eta_2^2}\in (0,1).
\end{equation}

$(2)$\; Solutions
$(\phi_0(\cdot;\eta_1,\eta_2),\psi_0(\cdot;\eta_1,\eta_2))$ given
by (\ref{eq:2.5}) and determined by $\eta_1=\eta_1(\sigma)$,
$\eta_2=\eta_2(\sigma)=\Lambda(\sigma)$, with
$\eta_1^2+\eta_2^2=4c\sigma$, have fundamental period $L$ and
satisfy   (\ref{eq:2.1}). Moreover, the mapping
$$
\sigma\in I(\sigma_0)\to \phi_0(\cdot;
\eta_1(\sigma),\eta_2(\sigma))\in
 H^n_{per}([0,L])
$$
is a smooth function (for all $n\geq 1$ integer).

$(3)$\; $I(\sigma_0)$ can be chosen as $(\frac{2\pi^2}{L^2},
+\infty)$.

$(4)$\; The mapping $\Lambda:I(\sigma_0)\to B(\eta_{2,0})$ is a
strictly decreasing function. Therefore, from (\ref{eq:2.13}),
$\sigma\to k(\sigma)$ is a strictly increasing function.
 \label{thm2.1}
\end{thm}

\noindent {\bf{Proof.}} The key of the proof is to apply the
implicit function theorem. In fact, consider the open set
$\Omega=\{(\eta,\sigma):\; \sigma>\frac{2\pi^2}{L^2},\;\eta\in
(0,\sqrt {2c\sigma})\;\}\subseteq \Bbb R^2$ and define
$\Psi:\Omega\to \Bbb R$ by
\begin{equation}\label{eq:2.14}
\Psi(\eta, \sigma)=\frac{4\sqrt
c}{\sqrt{4c\sigma-\eta^2}}\;K(k(\eta,\sigma))
\end{equation}
where $k^2(\eta,\sigma)=\frac{4c\sigma-2\eta^2}{4c\sigma-\eta^2}$.
By hypotheses $\Psi(\eta_{2,0},\sigma_0)=L$. Next, we show
$\partial_\eta\Psi(\eta,\sigma)<0$. In fact, it is immediate that
$$
\partial_\eta\Psi(\eta,\sigma)=\frac{4\sqrt
c\;\eta}{(4c\sigma-\eta^2)^{3/2}}\;K(k)+\frac{4\sqrt
c}{\sqrt{4c\sigma
-\eta^2}}\frac{dK}{dk}\frac{dk}{d\eta}.
$$
Next, from
$$
\frac{dk}{d\eta}=-\frac{4c\sigma\eta}{k(4c\sigma-\eta^2)^2},
$$
and the relations (see \cite{bf:bf})
\begin{equation}\label{eq:2.15}
\left \{
    \begin{array}{l}
\frac{dE}{dk}=\frac{E-K}{k},\qquad \frac{d^2E}{dk^2}=-\frac{1}{k}\frac{dK}{dk},\\
k{k'}^2\frac{d^2E}{dk^2}+{k'}^2\frac{dE}{dk}+kE=0,
\end{array}
   \right.
\end{equation}
with ${k'}^2=1-k^2$, and $E=E(k)$ being the complete elliptic
integral of second kind defined as
$$
E(k)=\int_0^1\;\sqrt{\frac{1-k^2t^2}{1-t^2}}\;dt,
$$
we have the following formal equivalences
\begin{eqnarray*}
 &\partial_\eta\Psi(\eta,\sigma)<0\Leftrightarrow
k(4c\sigma-\eta^2)\Big(E-k\frac{dE}{dk}\Big)<-4c\sigma k\frac{d^2E}{dk^2}\\
&\Leftrightarrow
k(4c\sigma-\eta^2)\Big(E-k\frac{dE}{dk}\Big)<\Big(\frac{dE}{dk}+\frac{k}{{k'}^2}
E \Big)(4c\sigma-\eta^2)(2-k^2)\\
&\Leftrightarrow 2{k'}^2 \frac{dE}{dk}+k E>0\Leftrightarrow
\frac{dE}{dk}-k\frac{d^2E}{dk^2}>0 \Leftrightarrow
\frac{dE}{dk}+\frac{dK}{dk}>0.
\end{eqnarray*}
So, since $E+K$ is a strictly increasing function, we obtain our
affirmation.

Therefore, there is a unique smooth function, $\Lambda$, defined
in a neighborhood $I(\sigma_0)$ of $\sigma_0$, such that
$\Psi(\Lambda(\sigma),\sigma)=L$ for every $\sigma\in
I(\sigma_0)$. So, we obtain (\ref{eq:2.12}). Moreover, since
$\sigma_0$ was chosen arbitrarily  in $\mathcal I=
(\frac{2\pi^2}{L^2}, +\infty)$, it follows from the uniqueness of
the function $\Lambda$ that it can be extended to $\mathcal I$.

Next, we show that $\Lambda$ is a strictly decreasing function. We
know that $\Psi(\Lambda(\sigma),\sigma)=L$ for every $\sigma\in
I(\sigma_0)$, then
\begin{eqnarray*}
\frac{d}{d\sigma}\Lambda(\sigma)=-\frac{\partial\Psi/{\partial
\sigma }}{\partial\Psi/{\partial\eta}}<0\Leftrightarrow
\partial\Psi/{\partial \sigma} <0.
\end{eqnarray*}
Thus, using the relation $\eta^2=(4c\sigma-\eta^2)(1-k^2)\equiv
(4c\sigma-\eta^2) {k'}^2$, we obtain the following formal
equivalences
$$
\frac{\partial\Psi}{\partial \sigma}<0 \Leftrightarrow
(4c\sigma-\eta^2)K>\frac{\eta^2}{k}\frac{dK}{dk}\Leftrightarrow
K>\frac{{k'}^2}{k} \frac{dK}{dk}.
$$
Then, since $\frac{dK}{dk}=(E-{k'}^2K)/(k{k'}^2)$, it follows that
$$
\frac{\partial\Psi}{\partial \sigma}<0 \Leftrightarrow k^2
K>E-{k'}^2K\Leftrightarrow (k^2+{k'}^2)K>E\Leftrightarrow K>E.
$$
This completes the proof of the Theorem. \hfill $\square$

\vspace{1.5em}

\noindent {\bf{Remark 4.1.}} In the case that the polynomial
$F(t)=-t^4+4c\sigma t^2+4cB_{\phi_0}$ has  a pure imaginary root and
the other two roots are real we can show the existence of  two
smooth curves of periodic solutions for (\ref{eq:2.2}) of cnoidal
type, more precisely we have
\begin{itemize}
\item $\omega\in (0,+\infty)\to b\;cn\Big(\sqrt{b^2-\omega}\; \xi
;k\Big)\in
 H^1_{per}([0,L])$
\item $\omega\in \Big(-\frac{4\pi^2}{L^2},0\Big)\to
\sqrt{a^2+2\omega}\;cn\Big(\sqrt{a^2+\omega}\; \xi ;k\Big)\in
 H^1_{per}([0,L])$,
\end{itemize}
where $a,b,k$ are smooth functions of $\omega$.

\vspace{1.5em}

The following result will be used in our stability theory.

\begin{cor}  Let $L$ and $c$ be arbitrarily fixed positive
numbers. Consider the smooth curve of dnoidal waves $\sigma\in
(\frac{2\pi^2}{L^2},\infty)\to \phi_0(\cdot;
\eta_1(\sigma),\eta_2(\sigma))$ determined by Theorem \ref{thm2.1}.
Then
$$
\frac{d}{d\sigma}\int_0^L \phi_0^2(\xi)\;d\xi>0.
$$
\label{cor2.2}
\end{cor}

\noindent {\bf{Proof.}} By (\ref{eq:2.5}), (\ref{eq:2.12}), and the
formula $\int_0^{K(k)} dn^2(x;k)\;dx=E(k)$ (see page 194 in
\cite{bf:bf}) it follows that
\begin{eqnarray*}
&\int_0^{L}\phi_0^2(\xi)\;d\xi= 2\eta_1\sqrt c\;
\int_0^{2K(k)}dn^2(x;k)\;dx=\frac{16cK}{L}\; \int_0^{K(k)}dn^2(x;k)\;dx\\
&=\frac{16c}{L}E(k)K(k).
\end{eqnarray*}
So, since
$k\to K(k) E(k)$ and $\sigma\to k(\sigma)$ are strictly increasing
functions we have that
$$
\frac{d}{d\sigma}\int_0^L \phi_0^2(\xi)\;d\xi=
\frac{16c}{L}\frac{d}{dk}[K(k) E(k)]\frac{dk}{d\sigma} >0.
$$
This finishes the Corollary. \hfill $\square$

\indent
\begin{center}
{4.2 \ \it Periodic Traveling Waves Solutions for Eq. (\ref{eq:eds})
}
\end{center}

In this subsection we show the existence of a branch of periodic
traveling waves solutions of (\ref{eq:eds}) for $\gamma$ close to
zero such that these solutions bifurcate the dnoidal waves solutions
found in Theorem \ref{eq:2.1}.

We start our analysis by studying  the periodic eigenvalue problem
considered on $[0,L]$,
\begin{equation}\label{eq:2.16}
\left \{
    \begin{array}{l}
 \mathcal L_{0}\chi\equiv (-\frac{d^2}{dx^2}+\sigma-\frac{3}{2c}\phi_0^2)\chi=\lambda\chi\\
 \chi(0)=\chi(L),\;\;\chi'(0)=\chi'(L),
 \end{array}
\right.
\end{equation}
where for $\sigma>2\pi^2/{L^2}$, $\phi_0$ is given by Theorem
\ref{eq:2.1} and satisfies \eqref{eq:2.2}.

\begin{thm}  The linear operator $\mathcal L_{0}$ defined in
(\ref{eq:2.16}) with domain $H^2_{per}([0,L])$ $\subseteq
L^2_{per}([0,L])$,  has its first three  eigenvalues simple with
zero being its second eigenvalue (with eigenfunction
$\frac{d}{dx}\phi_0$). Moreover, the remainder  of the spectrum is
constituted by a discrete set of eigenvalues which are double and
converging to infinity. \label{thm2.2}
\end{thm}

Theorem \ref{thm2.2} is a consequence of the Floquet theory
(Magnus\&Winkler \cite{mw:mw}). For convenience of the readers, we
will give some basic results of this theory. From the classical
theory of compact symmetric linear operator we have that problem
(\ref{eq:2.16}) determines a countable infinity set of eigenvalues
$\{\lambda_n|n=0,1,2,...\}$ with $
\lambda_0\leq\lambda_1\leq\lambda_2\leq\lambda_3\leq\lambda_4\leq...$,
where double eigenvalue is counted twice and $\lambda_n\to\infty$ as
$n\to\infty$. We shall denote by $\chi_n$ the eigenfunction
associated to the eigenvalue $\lambda_n$.  By the conditions $
\chi_n(0)=\chi_n(L),\;\;\chi'_n(0)=\chi'_n(L)$, $\chi_n$ can be
extended to the whole of  $(-\infty,\infty)$ as a continuously
differentiable function with period $L$.

We know that with the {\it periodic eigenvalue problem}
(\ref{eq:2.16}) there is an associated semi-periodic eigenvalue
problem in $[0,L]$, namely,
\begin{equation}\label{eq:2.17}
\left \{
    \begin{array}{l}
 \mathcal L_{0}\xi=\mu\xi\\
 \xi(0)=-\xi(L),\;\;\xi'(0)=-\xi'(L).
 \end{array}
\right.
\end{equation}
As in the periodic case, there is a countable infinity set of
eigenvalues $\{\mu_n|n=0,1,2,3,...\}$, with $
\mu_0\leq\mu_1\leq\mu_2\leq\mu_3\leq\mu_4\leq... $, where double
eigenvalue is counted twice and $\mu_n\to\infty$ as $n\to\infty$. We
shall denote by $\xi_n$ the eigenfunction associated to the
eigenvalue $\mu_n$.  So, we have that the equation
\begin{equation}\label{eq:2.18}
\mathcal L_{0}f=\gamma f
\end{equation}
has a solution of period $L$ if and only if $\gamma=\lambda_n$,
$n=0,1,2,\cdots$, while the only periodic solutions of period $2L$
are either those associated with $\gamma=\lambda_n$, but viewed on
$[0,2L]$, or those corresponding to  $\gamma=\mu_n$, but extended in
form $\xi_n(L+x)=\xi_n(L-x)$ for $0\leq x\leq L$, $n=0,1,2,\cdots$.
If all solutions of (\ref{eq:2.18}) are bounded we say that they are
{\it stable}; otherwise we say that they are {\it unstable}. From
the Oscillation Theorem of the Floquet theory (see \cite{mw:mw}) we
have that
\begin{equation}\label{eq:2.19}
\lambda_0<\mu_0\leq \mu_1<\lambda_1\leq \lambda_2<\mu_2\leq
\mu_3<\lambda_3\leq \lambda_4\cdot\cdot\cdot.
\end{equation}
The intervals $ (\lambda_0,\mu_0),
(\mu_1,\lambda_1),\cdot\cdot\cdot$, are called {\it intervals of
stability}. At the endpoints of these intervals the solutions of
(\ref{eq:2.18}) are unstable in general. This is always true for
$\gamma=\lambda_0$ ($\lambda_0$ is always simple).  The intervals,
$(-\infty,\lambda_0), (\mu_0,\mu_1)$, $(\lambda_1,\lambda_2),
(\mu_2,\mu_3),\cdots$, are called {\it intervals of
instability}\footnote{Here we omit any empty interval obtained from
a double eigenvalue.}. The interval $(-\infty,\lambda_0)$ of
instability will always be present. We note that {\it the absence of
an instability interval means that there is a value of $\gamma$ for
which all solutions of (\ref{eq:2.18}) have either period $L$ or
semi-period $L$, in other words, coexistence of solutions of
(\ref{eq:2.18}) with period $L$ or period $2L$ occurs for that value
of $\gamma$}. \vspace{1.5em}

\noindent {\bf{Proof of Theorem \ref{thm2.2}.}} From (\ref{eq:2.19})
we have that $\lambda_0<\lambda_1\leq \lambda_2$. Since $\mathcal
L_{0}\frac{d}{dx}\phi_0=0$ and $\frac{d}{dx}\phi_0$ has $2$ zeros in
$[0,L)$, it follows that $0$ is either $\lambda_1$ or $\lambda_2$.
We will show that $0=\lambda_1<\lambda_2$. We consider
$\Psi(x)\equiv\chi(\gamma x)$ with $\gamma^2=4c/{\eta_1^2}$. Then
from (\ref{eq:2.16}) and from the identity  $k^2 sn^2 x + dn^2 x=1$,
we obtain
\begin{equation}\label{eq:2.21}
\left \{
    \begin{array}{l}
 \frac{d^2}{dx^2}\Psi+[\rho-6k^2sn^2(x;k)]\Psi=0\\
\Psi(0)=\Psi(2K(k)),\;\;\Psi'(0)=\Psi'(2K(k)),
\end{array}
\right.
\end{equation}
where
\begin{equation}\label{eq:2.22}
\rho=\frac{4c(\lambda-\sigma)}{\eta_1^2}+6.
\end{equation}
Now, from Floquet theory, it follows that $(-\infty,\rho_0),
(\mu'_0,\mu'_1)$ and $(\rho_1,\rho_2)$ are the instability intervals
associated to this Lam\'e's equation, where for $i\geq 0$, $\mu'_i$
are the eigenvalues associated to the semi-periodic problem.
Therefore, $\rho_0,\rho_1,\rho_2$ are simple eigenvalues for
(\ref{eq:2.21}) and the other eigenvalues $\rho_3\leq
\rho_4<\rho_5\leq \rho_6<\cdot \cdot\cdot$ satisfy $\rho_3=\rho_4,
\rho_5= \rho_6,\cdot \cdot\cdot$, i.e., they are double eigenvalues.

It is easy to verify that the first three eigenvalues
$\rho_0,\rho_1,\rho_2$ and its corresponding eigenfunctions
$\Psi_0, \Psi_1, \Psi_2$ are given by the formulas
\begin{equation}\label{eq:2.23}
\left \{\begin{array}{l}
 \rho_0=2[1+k^2-\sqrt{1-k^2+k^4}],\;\;
 \Psi_0(x)=1-(1+k^2-\sqrt{1-k^2+k^4}\;)sn^2(x),\\
\rho_1=4+k^2,\;\;
\Psi_1(x)=sn x\; cn x,\\
\rho_2=2[1+k^2+\sqrt{1-k^2+k^4}],\;\;
\Psi_2(x)=1-(1+k^2+\sqrt{1-k^2+k^4}\;)sn^2(x).
\end{array}
\right.
\end{equation}
Next, $\Psi_0$ has no zeros in $[0,2K]$ and $\Psi_2$ has exactly 2
zeros in $[0,2K)$, then $\rho_0$ is the first eigenvalue to
(\ref{eq:2.21}). Since $\rho_0<\rho_1$ for every $k^2\in (0,1)$,
we obtain from (\ref{eq:2.22}) and (\ref{eq:2.6}) that
$$
4c\lambda_0=\eta_1^2(k^2-2-2\sqrt{1-k^2+k^4})<0\Leftrightarrow
\rho_0<\rho_1.
$$
Therefore $\lambda_0$ is the first negative eigenvalue to
$\mathcal L_{0}$ with eigenfunction $\chi_0(x)=\Psi_0(x/\gamma )$.
Similarly, since $\rho_1<\rho_2$ for every $k^2\in (0,1)$, we
obtain from (\ref{eq:2.22}) that
$$
4c\lambda_2=\eta_1^2(k^2-2+2\sqrt{1-k^2+k^4})>0\Leftrightarrow
\rho_1<\rho_2.
$$
Hence $\lambda_2$ is the third eigenvalue to $\mathcal L_{0}$ with
eigenfunction $\chi_2(x)=\Psi_2(x/\gamma )$. Finally, since
$\chi_1(x)=\Psi_1(x/\gamma)=\beta \frac{d}{dx}\phi_0(x)$ we finish
the proof. \hfill $\square$

\medskip

Next, we have our theorem of existence of solutions for
(\ref{eq:eds}). For $s\geq 0$, let $H^s_{per,e}([0,L])$ denote the
closed subspace of all even functions in $H^s_{per}([0,L])$.

\begin{thm} Let $L, \alpha,\beta, c>0$ and
$\sigma>2\pi^2/{L^2}$ be fixed numbers. Then there exist
$\gamma_1>0$ and a smooth branch
$$
\gamma\in (-\gamma_1,\gamma_1)\to
(\phi_{\gamma},\psi_\gamma)\in H^2_{per,e}([0,L])\times
H^1_{per,e}([0,L])
$$
of solutions for Eq. (\ref{eq:eds}). In particular, for $\gamma\to
0$, $(\phi_{\gamma},\psi_\gamma)$ converges to $(\phi_0,\psi_0)$
uniformly for $x\in [0,L]$, where $(\phi_0,\psi_0)$ is given by
Theorem \ref{thm2.1} and it is defined by (\ref{eq:2.5}). Moreover,
the mapping
$$
\gamma\in
(-\gamma_1,\gamma_1)\to (\frac{d}{d\sigma}
\phi_{\gamma},\frac{d}{d\sigma}\psi_\gamma)
$$
is continuous.
\label{thm2.3}
\end{thm}

\noindent {\bf{Proof.}} Without loss of generality, we take
$\alpha=1$ and $ \beta=1/2$. Let $X_e=H^2_{per,e}([0,L])\times
H^1_{per,e}([0,L])$ and define the map
$$
G:\Bbb R\times (0,+\infty)\times X_e\to
L^2_{per,e}([0,L])\times L^2_{per,e}([0,L])
$$
by
$$
G(\gamma, \lambda, \phi, \psi)=(-\phi''+\lambda \phi+\phi\psi,
-\gamma D\psi+c\psi+\frac12\phi^2).
$$
A calculation shows that the Fr\'echet derivative
$G_{(\phi,\psi)}=\partial
G(\gamma,\lambda,\phi,\psi)/{\partial(\phi,\psi)}$ exists and it is
defined as a map from $\Bbb R\times (0,+\infty)\times X_e$ to
$B(X_e;L^2_{per,e}([0,L])\times L^2_{per,e}([0,L]))$ by
\begin{eqnarray*}
G_{(\phi,\psi)}(\gamma, \lambda, \phi, \psi)=\begin{pmatrix}
  -\frac{d^2}{dx^2}+\lambda+\psi & \phi \\
  \phi & -\gamma D+c
\end{pmatrix}.
\end{eqnarray*}
>From Theorem \ref{thm2.1} it follows that for $\Phi_0=(\phi_0,
\psi_0)$, $G(0,\sigma, \Phi_0)=\vec 0^t$. Moreover, from Theorem
\ref{thm2.2} we have that $G_{(\phi,\psi)}(0, \sigma, \Phi_0)$ has a
kernel generated by ${\Phi_0'}^t$. Next, since $\Phi_0'\notin X_e$,
it follows that $G_{(\phi,\psi)}(0, \sigma, \Phi_0)$ is invertible.
Hence, since $G$ and $G_{(\phi,\psi)}$ are smooth maps on their
domains, the Implicit Function Theorem implies that there are
$\gamma_1>0$, $ \lambda_1\in (0, \sigma)$, and a smooth curve
$$
(\gamma,\lambda)\in (-\gamma_1,\gamma_1)\times
(\sigma-\lambda_1,\sigma+\lambda_1)\to (\phi_{\gamma,\lambda},
\psi_{\gamma,\lambda})\in X_e
$$
such that $G(\gamma,\lambda,\phi_{\gamma,\lambda},
\psi_{\gamma,\lambda})=0$. Then, for $\lambda=\sigma$ we obtain a
smooth branch $\gamma\in (-\gamma_1,\gamma_1)\to
(\phi_{\gamma,\sigma}, \psi_{\gamma,\sigma})\equiv (\phi_\gamma,
\psi_\gamma)$ of solutions of Eq. (\ref{eq:eds}) such that
$\gamma\in (-\gamma_1,\gamma_1)\to (\frac{d}{d\sigma}
\phi_{\gamma},\frac{d}{d\sigma}\psi_\gamma)$ is continuous. This
shows the Theorem. \hfill $\square$

\medskip

\noindent {\bf{Remark 4.2.}} Since $\phi_0$ is strictly positive on
$[0,L]$ and $\phi_\gamma\to \phi_0$, as $\gamma\to 0$, uniformly in
$[0,L]$, we have that, for $\gamma$ near zero, $\phi_\gamma(x)>0$
for $x\in \Bbb R$. Moreover, since the linear operator $-\gamma D+c$
is a strictly positive operator from $H^1_{per}([0,L])$ to
$L^2_{per}([0,L])$ for $\gamma$ negative, we have  that
$\psi_\gamma(x)<0$ for all $x\in \Bbb R$.

\setcounter{section}{4} \setcounter{equation}{0}
\begin{center}
\section{Stability of Periodic Traveling-Wave Solutions}
\end{center}

We begin this section defining the type of stability of our
interest. For any $c\in \Bbb R^+$ define the functions
$\Phi(\xi)=e^{ic\xi/2}\phi(\xi)$ and $\Psi(\xi)=\psi(\xi)$, where
$(\phi,\psi)$ is a solution of (\ref{eq:eds}). Then we say that the
orbit generated by $(\Phi,\Psi)$, namely,
$$
\Omega_{(\Phi,\Psi)}=\{(e^{i\theta}\Phi(\cdot+x_0),\Psi(\cdot+x_0))\;:
(\theta,x_0)\in [0,2\pi)\times \Bbb R\},
$$
is stable in $H^1_{per}([0,L])\times H^{\frac12}_{per}([0,L])$ by
the flow generated by Eq. (\ref{eq:lsws}), if for every $\epsilon
>0$ there exists $\delta (\epsilon)
>0$ such that for $(u_0, v_0)$ satisfying  $\|u_0-\Phi\|_1<\delta $
and $\|v_0-\Psi\|_\frac12<\delta $, we have that $(u,v)$ solution of
(\ref{eq:lsws}) with $(u(0),v(0))=(u_0,v_0)$, satisfies that
$(u,v)\in C(\mathbb R;H^1_{per}([0,L]))\times C(\mathbb R;
H^{\frac12}_{per}([0,L]))$ and
\begin{equation}\label{eq:3.1}
\inf\limits_{x_0\in\Bbb R,\theta\in [0,2\pi)}
\|e^{i\theta}u(\cdot+x_0,t)-\Phi\|_1 <\epsilon,\;\;
\inf\limits_{x_0\in\Bbb R}\;\|v(\cdot+x_0,t)-\Psi\|_{\frac12}
<\epsilon,
\end{equation}
for all $t\in \mathbb R$.

\vspace{1.5em}

The main result to be proved in this section is that the  periodic
traveling waves solutions of (\ref{eq:lsws}) determined by Theorem
\ref{thm2.3} are stable for $\sigma>2\pi^2/{L^2} $ and $\gamma$
negative close to 0.

\medskip

\begin{thm} Let $L, \alpha,\beta, c>0$ and
$\sigma>2\pi^2/{L^2}$ be fixed numbers. We consider the smooth curve
of periodic traveling waves solutions for (\ref{eq:eds}), $\gamma\to
(\phi_\gamma,\psi_\gamma)$, determined by Theorem \ref{thm2.3}. Then
there exists $\gamma_0>0$ such that for each $\gamma\in (-\gamma_0,
0)$, the orbit generated by $(\Phi_\gamma(\xi), \Psi_\gamma(\xi))$
with
$$
\Phi_\gamma(\xi)=e^{ic\xi/2}\phi_\gamma(\xi)\;\;\text{and}\;\;
\Psi_\gamma(\xi)=\psi_\gamma(\xi),
$$
is orbitally stable in $H^1_{per}([0,L])\times
H^{\frac12}_{per}([0,L])$.
 \label{thm3.1}
\end{thm}

The proof of Theorem \ref{thm3.1} is based on the ideas developed by
Benjamin (\cite{b:b}) and  Weinstein (\cite{w1:w1}) which give us an
easy form of manipulating with the required  spectral information
and the positivity property of the quantity
$\frac{d}{d\sigma}\int\phi^2_\gamma(x)dx$, which are basic in our
stability theory. We do not use the abstract stability theory of
Grillakis {\it et al.} basically by these circumstance. So, consider
$(\phi_\gamma, \psi_\gamma)$ a solution of (\ref{eq:eds}) obtained
in Theorem \ref{eq:2.3}. For $(u_0,v_0)\in H^{1}_{per}([0,L])\times
H^{\frac12}_{per}([0,L])$ and $(u,v)$ the global solution to
(\ref{eq:lsws}) corresponding to these initial data given by Theorem
\ref{theog}, we define for $t\geq 0$ and $\sigma>2\pi^2/{L^2}$
\begin{eqnarray}\label{eq:3.2}
\Omega_t(x_0,\theta)= \| e^{i\theta}(T_c u)'(\cdot+x_0,t)-
\phi_\gamma'\|^2
  + \sigma \|e^{i\theta}(T_cu)(\cdot+x_0,t)- \phi_\gamma\|^2
\end{eqnarray}
where we denote by $T_c$ the bounded linear operator  defined by
$$
(T_c u)(x,t)=e^{-ic(x-ct)/2}u(x,t).
$$
Then, the deviation of the solution $u(t)$ from the orbit
generated by $\Phi$ is measured by
\begin{eqnarray}\label{eq:3.3}
\rho_\sigma(u(\cdot,t),
\phi_\gamma)^2\equiv&\;\;\inf\limits_{x_0\in [0,L],\theta\in
[0,2\pi]}\;\Omega_t(x_0,\theta).
\end{eqnarray}
Hence, from (\ref{eq:3.3}) we have that the $ \inf
\Omega_{t}(x_0,\theta)$ is attained in $(\theta, x_0)=(\theta(t),
x_0(t))$.

\vspace{1.0em}

\noindent {\bf{Proof of Theorem \ref{thm3.1}.}}  Consider the
perturbation of the periodic traveling wave
$(\phi_\gamma,\psi_\gamma)$
\begin{equation}\label{eq:3.4}
\left \{\begin{array}{l}
\xi (x,t)=e^{i\theta}(T_c u)(x+x_0,t)-\phi_{\gamma}(x) \\
\eta (x,t)= v(x+x_0,t)-\psi_{\gamma}(x).
\end{array}
\right.
\end{equation}
Hence, by the property of minimum of $(\theta, x_0)=(\theta(t),
x_0(t))$, we obtain from (\ref{eq:3.4}) that $p(x,t)=\text{Re}(\xi
(x,t))$ and $q(x,t)=\text{Im}(\xi (x,t))$ satisfy the compatibility
relations
\begin{equation}\label{eq:3.5}
\left \{\begin{array}{l}
\int_0^L q(x,t)\phi_\gamma(x)\psi_\gamma(x)\;dx=0\\
\int_0^L p(x,t)(\phi_\gamma(x)\psi_\gamma(x))'\;dx=0.
\end{array}
\right.
\end{equation}

Now we take the continuous functional $L$ defined on
$H^1_{per}([0,L])\times H^{\frac12}_{per}([0,L])$ by
$$
L(u,v)= E(u,v)+ c\; G(u,v) + \omega\; H(u,v),
$$
where $E, G, H$ are defined by (\ref{eq:law2}). Then, from
(\ref{eq:3.4}) and (\ref{eq:eds}), we have
\begin{equation}\label{eq:3.6}
\begin{array}{l}
\Delta L(t):= L(u(t),v(t))- L(\Phi_{\gamma}, \Psi_{\gamma})=
L(\Phi_{\gamma}+e^{icx/2}\xi, \psi_{\gamma}+ \eta)- L(\Phi_{\gamma}, \psi_{\gamma})\\
\\
= <\mathcal L_\gamma p,p> + <\mathcal L_\gamma^+ q,q> +
\frac{\alpha}{2\beta}\int_{0}^L \Big [{\mathcal
K}_{\gamma}^{1/2}\eta+ 2\beta {\mathcal
K}_{\gamma}^{-1/2}(\phi_{\gamma}p)+\beta {\mathcal
K}_{\gamma}^{-1/2}(
p^2+q^2)\Big]^2dx\\
\\
-\frac{\alpha \beta}{2}\int_{0}^L \Big[|{\mathcal
K}_{\gamma}^{-1/2} (p^2+q^2)|^2+4{\mathcal
K}_{\gamma}^{-1/2}(\phi_{\gamma}p){\mathcal K}_{\gamma}^{-1/2}
(p^2+q^2)\Big]dx,
\end{array}
\end{equation}
where, for $\gamma<0$ we define $\mathcal K_\gamma^{-1}$  as
$$
\widehat {\mathcal K_\gamma^{-1}f}(k)=\frac{1}{-\gamma
|k|+c}\widehat f(k)\qquad\text{for}\; k\in\Bbb Z,
$$
which is the inverse operator of $\mathcal K_\gamma :
H^s_{per}([0,L])\to H^{s-1}_{per}([0,L])$  defined by $\mathcal
K_\gamma= -\gamma D+c$. The operator $\mathcal L_\gamma$ is
\begin{equation}\label{eq:3.7}
\mathcal L_\gamma= -\frac{d^2}{d\xi^2} +\sigma
+\alpha\psi_\gamma-2\alpha\beta \phi_\gamma\circ\mathcal
K_\gamma^{-1}\circ\phi_\gamma,
\end{equation}
with $\phi_\gamma\circ\mathcal K_\gamma^{-1}\circ\phi_\gamma$ given
by $[\phi_\gamma\circ\mathcal K_\gamma^{-1}\circ\phi_\gamma]
(f)=\phi_\gamma\mathcal K_\gamma^{-1}(\phi_\gamma f)$. Here
$\mathcal L_\gamma^+$ is defined by
\begin{equation}\label{eq:3.8}
\mathcal L_\gamma^+ = -\frac{d^2}{d \xi^2}+\sigma  + \alpha
\psi_{\gamma}
\end{equation}
and ${\mathcal K}_{\gamma}^{1/2}$, ${\mathcal K}_{\gamma}^{-1/2}$
are the positive roots of ${\mathcal K}_{\gamma}$ and ${\mathcal
K}_{\gamma}^{-1}$ respectively.

Now, we need to find a lower bound for $\Delta {L}(t)$. The first
step will be to obtain a suitable lower bound of the last term on
the right-hand side of (\ref{eq:3.6}). In fact, since ${\mathcal
K}_{\gamma}^{-1/2}$ is a bounded operator on $L^2_{per}([0,L])$,
$\phi_\gamma$ is uniformly bounded, and from the continuous
embedding of $H^1_{per}([0,L])$ in $L^4_{per}([0,L])$ and in
$L^{\infty}([0,L])$, we have that
\begin{equation}\label{eq:3.9}
-\frac{\alpha \beta}{2}\int_{0}^L \Big [|{\mathcal
K}_{\gamma}^{-1/2} (p^2+q^2)|^2+4{\mathcal K}_{\gamma}^{-1/2}
(\phi_{\gamma}p){\mathcal K}_{\gamma}^{-1/2}(p^2+q^2)\Big ]dx \geq
-C_1\|\xi \|_1^3-C_2\|\xi \|_1^4
\end{equation}
where $C_1$ and $C_2$ are positive constants.

\vspace{1.5em}

 The estimates for $\langle\mathcal L_\gamma p,p\rangle$ and
$\langle\mathcal L_\gamma^+ q,q\rangle$ will be obtained from the
following theorem.

\begin{thm} Let $L, \alpha,\beta, c>0$ and
$\sigma>2\pi^2/{L^2}$ be fixed numbers. Then, there exists
$\gamma_2>0$ such that, if $\gamma\in (-\gamma_2,0)$, the
self-adjoint operators $\mathcal L_\gamma$ and $\mathcal L^+_\gamma$
defined in (\ref{eq:3.7}) and (\ref{eq:3.8}), respectively, with
domain $H^2_{per}([0,L])$ have the following properties:

$(1)$  $\mathcal L_\gamma$ has a simple negative eigenvalue
$\lambda_\gamma$ with eigenfunction $\varphi_\gamma$ and
$\int_{0}^L\phi_\gamma\varphi_\gamma dx \neq 0$.

$(2)$ $\mathcal L_\gamma$ has a simple eigenvalue at zero with
eigenfunction $\frac{d}{dx}\phi_\gamma$.

$(3)$ There is $\eta_\gamma>0$ such that for $\beta_\gamma\in
\Sigma(\mathcal L_\gamma)-\{\lambda_\gamma,0\}$, we have that
$\beta_\gamma>\eta_\gamma$.

$(4)$  $\mathcal L^+_\gamma$ is a non-negative operator which has
zero as its first eigenvalue with eigenfunction $\phi_\gamma$. The
remainder of the spectrum is constituted by a discrete set of
eigenvalues.
 \label{thm3.2}
\end{thm}

\noindent {\bf{Proof.}} From (\ref{eq:eds}) it follows that
$\mathcal L_\gamma\phi_\gamma=2\phi_\gamma\psi_\gamma$ and so, from
Remark 4.2, we have that for $\gamma<0$,
 $\langle\mathcal
L_\gamma \phi_\gamma,\phi_\gamma\rangle=2\int_{\Bbb
R}\phi^2_\gamma\psi_\gamma dx<0$. Therefore $\mathcal L_\gamma$ has
a negative eigenvalue. Moreover, we have that $\mathcal L_\gamma
\frac{d}{dx}\phi_\gamma=0$. Next, for $f\in H^1_{per}([0,L])$ and
$\|f\|=1$, we have
\begin{equation}\label{eq:3.10}
\begin{array}{l}
 \langle\mathcal L_\gamma f,f\rangle=\langle\mathcal L_0
 f,f\rangle-\frac{\gamma}{c}\alpha^2\langle\phi_0 f, D\mathcal
 K^{-1}_{\gamma}(\phi_0f)\rangle+\alpha \int_{0}^L(\psi_\gamma-\psi_0)f^2dx\\
\\
+\alpha^2 \int_{0}^L[\phi_0f \mathcal
 K^{-1}_{\gamma}(\phi_0f)-\phi_\gamma f \mathcal
 K^{-1}_{\gamma}(\phi_\gamma f)]dx\\
 \\
 \geq \langle\mathcal L_0
 f,f\rangle+ \alpha \int_{0}^L(\psi_\gamma-\psi_0)f^2dx+\alpha^2 \int_{0}^L[\phi_0f
 \mathcal
 K^{-1}_{\gamma}(\phi_0f)-\phi_\gamma f \mathcal
 K^{-1}_{\gamma}(\phi_\gamma f)]dx,
\end{array}
\end{equation}
where the last inequality is due to that $\gamma<0$ and $D\mathcal
 K^{-1}_{\gamma}$ is a positive operator. So, since
\begin{equation}\label{eq:3.11}
\begin{array}{l}
\Big|\int_{0}^L(\psi_\gamma-\psi_0)f^2dx\Big|\leq
 \|\psi_\gamma-\psi_0\|_{\infty}\\
 \\
\Big|\int_{0}^L[\phi_0f \mathcal
 K^{-1}_{\gamma}(\phi_0f)-\phi_\gamma f \mathcal
 K^{-1}_{\gamma}(\phi_\gamma f)]dx\Big|\leq
 (\|\phi_\gamma\|+\|\phi_0\|)\|\phi_\gamma-\phi_0\|_{\infty},
 \end{array}
\end{equation}
we have from Theorem \ref{thm2.3} that for $\gamma$ near $0^-$ and
$\epsilon$ small, $\langle\mathcal L_\gamma
f,f\rangle\geq\langle\mathcal L_0
 f,f\rangle-\epsilon$. Hence, for $f\perp \chi_0$ and $f\perp \frac{d}{dx}\phi_0$, where $\mathcal
L_0 \chi_0=\lambda_0 \chi_0$ with $\lambda_0<0$, we have from the
spectral structure of $\mathcal L_0$ (Theorem \ref{thm2.2}) that $
\langle\mathcal L_\gamma f,f\rangle\geq \eta_\gamma>0$. Therefore,
from min-max principle (\cite{rs:rs}) we obtain the desired spectral
structure for $\mathcal L_\gamma $. Moreover, let $\varphi_\gamma$
be such that $\mathcal
 L_\gamma
\varphi_\gamma=\lambda_\gamma \varphi_\gamma$ with
$\lambda_\gamma<0$. Therefore, if $\phi_\gamma \perp\varphi_\gamma$
then from the spectral structure of $\mathcal L_\gamma$ we must have
that $\langle\mathcal L_\gamma \phi_\gamma,\phi_\gamma\rangle\geq
0$. But we know that $\langle\mathcal L_\gamma
\phi_\gamma,\phi_\gamma\rangle <0$. Hence, $\langle\phi_\gamma,
\varphi_\gamma\rangle\neq 0$. Finally, since $\mathcal
L_\gamma^+\phi_\gamma=0$ with $\phi_\gamma>0$, it follows that zero
is simple and it is the first eigenvalue. The remainder of the
spectrum is discrete. \hfill $\square$

\vspace{1.5em}

\begin{thm}Consider $\gamma<0$ close to zero such that Theorem
\ref{thm3.2} is true. Then
\begin{eqnarray*}
(a)&\inf\;\{\langle\mathcal L_\gamma f,f\rangle\;:\;\|f\|=1,
\langle f,\phi_\gamma\rangle=0,\;\}\equiv \beta_0=0.\\
 (b)&\inf\;\{\langle\mathcal
L_\gamma f,f\rangle\;:\;\|f\|=1, \langle f,\phi_\gamma\rangle=0,
\langle f,(\phi_\gamma\psi_\gamma)'\rangle=0\;\}\equiv \beta>0.
\end{eqnarray*}
\label{thm3.3}
\end{thm}

\noindent {\bf{Proof.}} \textit{Part (a).} Since $\mathcal L_\gamma
\frac{d}{dx}\phi_\gamma=0$ and $\langle\frac{d}{dx}\phi_\gamma,
\phi_\gamma\rangle=0$ then $\beta_0\leq 0$. Next we will show that
$\beta_0\geq 0$ by using Lemma E.1 in Weinstein \cite{w2:w2}. So, we
shall show initially that the infimum  is attained. Let
$\{\psi_j\}\subseteq H^1_{per}([0,L])$ with $\|\psi_j\|=1$,
$\langle\psi_j,\phi_\gamma\rangle=0$ and $\lim_{j\to
\infty}\langle\mathcal L_\gamma \psi_j,\psi_j\rangle=\beta_0$. Then
there is a subsequence of $\{\psi_j\}$, which we denote again by
$\{\psi_j\}$, such that $\psi_j \rightharpoonup \psi$ weakly in
$H^1_{per}([0,L])$, so $\psi_j \to \psi$ in $L^2_{per}([0,L])$.
Hence $\|\psi\|=1$ and $\langle\psi,\phi_\gamma\rangle=0$. Since
$\|\psi'\|^2\leq \liminf \|\psi'\|^2$ and $\mathcal
K^{-1}_\gamma(\phi_\gamma \psi_j)\to \mathcal
K^{-1}_\gamma(\phi_\gamma \psi)$ in $L^2_{per}([0,L])$, we have $
\beta_0\leq \langle\mathcal L_\gamma \psi,\psi\rangle\leq \liminf
\langle\mathcal L_\gamma \psi_j,\psi_j\rangle=\beta_0$. Next we show
that $\langle\mathcal L^{-1}_\gamma
\phi_\gamma,\phi_\gamma\rangle\leq 0$. From (\ref{eq:eds}) and
Theorem \ref{thm2.3} we obtain  for $\chi_\gamma=-\frac{d}{d\sigma}
\phi_\gamma$ that $\mathcal L_\gamma \chi_\gamma =\phi_\gamma$.
Moreover, from Corollary \ref{cor2.2} it follows that
$\langle-\frac{d}{d\sigma} \phi_0,\phi_0\rangle<0$ and so for
$\gamma$ small enough $\langle-\frac{d}{d\sigma}
\phi_\gamma,\phi_\gamma\rangle<0$. Hence from \cite{w2:w2} we obtain
that $\beta\geq 0$. This shows part (a) of the Theorem.

\textit{Part (b).} From  (a) we have that $\beta\geq 0$. Suppose
$\beta=0$. Then  following a similar analysis to that used in part
(a) above, we have that the infimum  define in (b) is attained at an
admissible function $\zeta$. So, from Lagrange's multiplier theory,
there are $\lambda, \theta, \eta$ such that
\begin{equation}\label{eq:3.12}
\mathcal L_\gamma \zeta=\lambda \zeta+\theta \phi_{\gamma}+ \eta
(\phi_{\gamma}\psi_{\gamma})'.
\end{equation}
Using (\ref{eq:3.12}) and $\langle\mathcal L_\gamma
\zeta,\zeta\rangle=0$ we obtain that $\lambda =0$. Taking the inner
product of (\ref{eq:3.12}) with $\phi'_{\gamma}$, we have from
$\mathcal L_\gamma \phi'_{\gamma}=0$ that
\begin{equation}\label{eq:3.13}
0=\eta \int_{0}^L \phi'_{\gamma}(\phi_{\gamma}\psi_{\gamma})'dx,
\end{equation}
but the integral in (\ref{eq:3.13}) converges to
$$
\int_{0}^L \phi'_0(\phi_0\psi_0)'dx=\frac{-3\beta}{c}\int_{0}^L
\phi_0^2(\phi_0')^2dx <0
$$
as $\gamma \to 0$. Then, from (\ref{eq:3.13}), we obtain $\eta =0$
and therefore $ \mathcal L_\gamma \zeta=\theta \phi_{\gamma}$. So,
since $\mathcal L_\gamma (-\frac{d}{d\sigma} \phi_\gamma)
=\phi_\gamma$, we obtain $
0=\langle\zeta,\phi_{\gamma}\rangle=\theta \langle\phi_{\gamma},
-\frac{d}{d\sigma} \phi_\gamma\rangle$. Therefore $\theta =0$ and
$\mathcal L_\gamma \zeta=0$. Then $\zeta=\nu \phi_{\gamma}'$ for
some $\nu \neq 0$, which is a contradiction. Thus $\beta>0$ and the
proof of the Theorem is completed. \hfill $\square$

\begin{thm}\label{thm3.4}Consider $\gamma<0$ close to zero such that Theorem
\ref{thm3.2} is true. If $\mathcal L^+_\gamma$ is defined as in
(\ref{eq:3.8}) then
\begin{eqnarray*}
\inf\;\{\langle\mathcal L^+_\gamma f,f\rangle\;:\;\|f\|=1, \langle
f,\phi_\gamma\psi_\gamma\rangle=0,\;\}\equiv \mu>0.
\end{eqnarray*}
\end{thm}

\noindent {\bf{Proof.}} From Theorem \ref{thm3.2} we have that
$\mathcal L_\gamma^+$ is a non-negative operator and so $\mu\geq 0$.
Suppose $\mu=0$. Then, by following the ideas of the proof of
Theorem \ref{thm3.3}, we have that the minimum is attained at an
admissible function $g^*\neq 0$ and there is $(\lambda , \theta )
\in \Bbb R^2$ such that
\begin{equation}\label{eq:3.14}
\mathcal L_\gamma^+ g^*=\lambda g^* +\theta
\phi_{\gamma}\psi_{\gamma}.
\end{equation}
Thus, it follows that $\lambda =0$. Now, taking the inner product of
(\ref{eq:3.14}) with $\phi_{\gamma}$ it is deduced that $
0=\langle\mathcal L_\gamma^+\phi_{\gamma},g^*\rangle=\theta
\int_{0}^L \phi_{\gamma}^2\psi_{\gamma}dx$, and therefore $\theta
=0$. Then, since zero is a simple eigenvalue for $\mathcal
L_\gamma^+$ it follows that $g^*=\nu \phi_{\gamma}$ for some $\nu
\neq 0$, which is a contradiction. This completes the proof. \hfill
$\square$

\vspace{1.5em}

Next we prove Theorem \ref{thm3.1} by returning to (\ref{eq:3.6}).
Our task is to estimate the terms $\langle\mathcal L_\gamma
p,p\rangle$ and $\langle\mathcal L_\gamma^+ q,q\rangle$ where $p$
and $q$ satisfy (\ref{eq:3.5}). From Theorem \ref{thm3.4} and the
definition of $\mathcal L_\gamma^+$, we have that there is $C_1>0$
such that
\begin{equation}\label{eq:3.15}
\langle\mathcal L_\gamma^+ q,q\rangle\geq C_1\|q\|_1^2.
\end{equation}
Now we estimate $\langle\mathcal L_\gamma p,p\rangle$. Suppose
without loss of generality that $\|\phi_\gamma\|=1$. We write
$p_{_{\perp}}= p-p_{_{||}}$, where $p_{_{||}}=\langle
p,\phi_\gamma\rangle\phi_\gamma$. Then, from (\ref{eq:3.5}) and the
positivity of the operator ${\mathcal K}_{\gamma}^{-1}$ it follows
that $\langle p_{_{\perp}},(\phi_{\gamma}\psi_{\gamma})'\rangle=0$.
Therefore from Theorem \ref{thm3.3}, it follows $ \langle\mathcal
L_\gamma p_{_{\perp}},p_{_{\perp}}\rangle\geq D\|p_{_{\perp}}\|^2$.
Now we suppose that $ \|u_0\|=\|\phi_{\gamma}\|=1$. Since
$\|u(t)\|^2=1$ for all $t$, we have that $\langle
p,\phi_\gamma\rangle=-\|\xi\|^2/2$. So, $\langle\mathcal L_\gamma
p_{_{\perp}},p_{_{\perp}}\rangle\geq
\beta_0\|p\|^2-\beta_1\|\xi\|_1^4$. Since $\langle\mathcal L_\gamma
\phi_{\gamma},\phi_{\gamma}\rangle\;<0$ it follows that
$\langle\mathcal L_\gamma p_{_{||}},p_{_{||}}\rangle\geq
-\beta_3\|\xi\|_1^4$. Moreover, Cauchy-Schwarz inequality implies
$\langle\mathcal L_\gamma p_{_{||}},p_{_{\perp}}\rangle\geq
-\beta_4\|\xi\|_1^3$. Therefore we conclude from the specific form
of $\mathcal L_\gamma$ that
\begin{equation}\label{eq:3.16}
\langle\mathcal L_\gamma p ,p\rangle\geq
D_1\|p\|_{1,\sigma}^2-D_2\|\xi\|_{1,\sigma}^3-D_3\|\xi\|_{1,\sigma}^4,
\end{equation}
with $D_i>0$ and $\|f\|_{1,\sigma}^2=\|f'\|^2+\sigma\|f\|^2$.

Next, by collecting the results in (\ref{eq:3.9}), (\ref{eq:3.15})
and  (\ref{eq:3.16}) and substituting them in (\ref{eq:3.6}), we
obtain
\begin{equation}\label{eq:3.17}
\Delta L(t) \geq
d_1\|\xi\|_{1,\sigma}^2-d_2\|\xi\|_{1,\sigma}^3-d_3\|\xi\|_{1,\sigma}^4,
\end{equation}
where $d_i>0$. Therefore, from standard arguments, for any
$\epsilon>0$, there exists $\delta(\epsilon)>0$ such that, if
$\|u_0-\Phi_\gamma\|_{1,\sigma}<\delta(\epsilon) $ and
$\|v_0-\Psi_\gamma\|_{\frac12}<\delta(\epsilon) $, then
\begin{equation}\label{eq:3.18}
\rho_\sigma(u(t), \phi_\gamma)^2=\|\xi(t)\|^2_{1,\sigma}<\epsilon
\end{equation}
for $t\in [0,\infty)$, and so we obtain the first inequality in
(\ref{eq:3.1}).

Now, it follows from (\ref{eq:3.6}) and from the above analysis of
$\xi$ that
$$
\epsilon\geq \frac{\alpha}{2\beta}\int_{\Bbb R} \Big [{\mathcal
K}_{\gamma}^{1/2}\eta+ 2\beta {\mathcal
K}_{\gamma}^{-1/2}(\phi_{\gamma}p)+\beta {\mathcal
K}_{\gamma}^{-1/2}( p^2+q^2)\Big]^2dx.
$$
Thus, from (\ref{eq:3.18}) and the equivalence of the norms
$\|{\mathcal K}_{\gamma}^{1/2}\eta\|$ and $\|\eta\|_{\frac12}$, we
obtain (\ref{eq:3.1}). This proves that $(\Phi_\gamma,\Psi_\gamma)$
is stable relative to small perturbation which preserves the
$L^2_{per}([0,L])$ norm of $\Phi_\gamma$. The general case  follows
from that $\gamma\in (-\gamma_1,\gamma_1)\to
(\phi_{\gamma},\psi_\gamma)$ is a smooth branch of solutions for Eq.
(\ref{eq:eds}). \hfill $\square$


\indent

\end{document}